\documentclass[a4paper,twoside,10pt]{amsart}

\usepackage{amssymb}
\usepackage{amsmath}
\usepackage{epsfig}

\usepackage{eepic}

\usepackage[latin1]{inputenc} 
\usepackage[T1]{fontenc}

\usepackage[active]{srcltx}

\usepackage{eepic}

\newtheorem{The}{Theorem}[section]
\newtheorem{Lem}[The]{Lemma}
\newtheorem{Cor}[The]{Corollary}
\newtheorem{Pro}[The]{Proposition}

\theoremstyle{definition}
\newtheorem{definition}[The]{Definition}
\newtheorem{Exam}[The]{Example}
\newtheorem{Not}[The]{Notation}

 \theoremstyle{remark}
\newtheorem{remark}[The]{Remark}

\def\a{\alpha}
\def\A{\mathbf A}

\def\Bar{B_{\mathrm {ar}}}
 
\def\cc{{\mathcal{C}}}

\def\f{\phi}

\def\J{{\mathcal{J}}}

\def\kvar{{K_0 (\mathrm{Var}_k)}}
\def\kchow{{K_0 (\mathrm{CHMot}_k)}}
\def\kmot{{K_0^{\mathrm{mot}}(\mathrm{Var}_k)}}
\def\kmotq {{K_0^{\mathrm{mot}}(\mathrm{Var}_k) \otimes\mathbf Q}}
\def\kfield{{K_0(\mathrm{Field}_k)}}
\def\kpff{{K_0(\mathrm{PFF}_k)}}

\def\l{\lambda}
\def\L{\mathbf{L}}

\def\Newton{\mathcal N}

\def\ord{\mbox{\rm ord}}
\def\span{\mbox{\rm span}}

\def\O{\mbox{\rm orb}}
\def\p{\pi}

\def\Pgeom{P_{\mathrm {geom}}}
\def\Par{P_{\mathrm {ar}}}
\def\Q{\mathbf Q}
\def\R{\mathbf R}
\def\r{\rho}

\def\s{\sigma}
\def\t{\tau}

\def\y{\wedge}

\def\Z{\mathbf Z}

\def\limproj{\mathop{\oalign{lim\cr\hidewidth$\longleftarrow$\hidewidth\cr}}}

\begin{document}

\title{Arithmetic Motivic Poincaré series of toric varieties}

\author{H. Cobo Pablos}

 \address{  Department of Mathematics,
     University of Leuven,
     Celestijnenlaan 200B,
B-3001 Leuven-Heverlee, Belgium }

 \email{Helena.Cobo@wis.kuleuven.ac.be}

\author{P.D. Gonz{á}lez P{é}rez}

\address{Instituto de Ciencias Matemáticas-CSIC-UAM-UC3M-UCM. Depto. \'Algebra. Facultad de Ciencias Matemáticas. Universidad Complutense de Madrid.
Plaza de las Ciencias 3. 28040. Madrid. Spain}

\email{pgonzalez@mat.ucm.es}

\thanks{ Helena Cobo Pablos is supported by FWO-Flanders project G031806N.
Pedro D. Gonz\'alez P\'erez is supported by {\em Programa Ram\'on y Cajal}
of {\em Ministerio de Ciencia e Innovación} (MCI), Spain.
Both authors are supported by MCI grant MTM2010-21740-C02-01.}



\begin{abstract}
     The \textit{arithmetic motivic Poincar\'e series} of a variety
$V$  defined over a field of characteristic zero, is an invariant of singularities
which was introduced by Denef and Loeser by analogy with the  Serre-Oesterlé series in arithmetic geometry.
They proved that this motivic series has a rational form which specializes to
the Serre-Oesterlé series when  $V$ is defined over the integers.
This invariant, which is known explicitly for a few classes of singularities, remains quite mysterious.
In this paper we study this motivic series when $V$ is an affine toric variety.
We obtain a formula for the rational form of this series in  terms of the Newton polyhedra of the  \textit{ideals of
sums of combinations} associated to the minimal system of generators of the semigroup of the toric variety. In particular, we
deduce explicitly a finite set of candidate poles for this invariant.
\end{abstract}

\maketitle

\section*{Introduction}

Let $S$ denote an irreducible and reduced algebraic variety defined over a field $k$ of characteristic zero.
The set $H(S)$ of
formal arcs of the form $ \mbox{Spec } k[[t]] \rightarrow S$ can
be given the structure of scheme over $k$ (not necessarily of
finite type). If $0 \in S$ we denote by $H(S)_0$ the subscheme of the arc space
consisting on arcs in $H(S)$ with origin at $0$.  The set $H_m(S)$ of $m$-jets of $S$,  of the form $
\mbox{Spec } k[t] / (t^{m+1}) \rightarrow S$, has the structure
of algebraic variety over $k$.
By a theorem of Greenberg, the
image of the space of arcs $H(S)$ by the natural morphism of
schemes $j_m : H(S) \rightarrow H_m(S)$ which maps any arc to
its $m$-jet,
  is a constructible subset of  $ H_m(S)$.

It follows from this that $j_m (H (S))$ defines a class $[j_m (H (S))]$ in the Grothendieck
ring of varieties $\kvar$ and also  a class
$\chi_c ([H_m (S)] ) \in \kchow$ in the Grothendieck ring of Chow motives, where  $\chi_c:\kvar\rightarrow\kchow$
is the unique ring homomorphism,
which maps the class of a smooth projective variety to its Chow motive (see \cite{GS,GN}).
 We denote by $\kmot$ the image of $\kvar$ by the homomorphism $\chi_c$. We use the
same symbol $\L$ to denote the class $[\A^1_k] \in \kvar$ and the class $\chi_c ([\A^1_k]) \in \kmot$.

Denef and Loeser have defined various notions of motivic Poincaré series,
 motivated by some generating series in arithmetic geometry (see \cite{DL-Dwork}).
Assume for simplicity that the variety $S$ is defined over the integers. We denote by $p$ a prime number and by $\Z_p$ the
$p$-adic integers. For every positive integer $m$, the symbol $N_{p,m}(S)$ denotes the number of
rational points of $S$ over  $\Z/p^{m+1}\Z$ which can be lifted to
rational points of $S$ over $\Z_p$ by the projection induced by the natural map   $\Z_p\rightarrow\Z/p^{m+1}\Z$.
The {\em Serre-Oesterlé} series of $S$ at the prime $p$ is
\[ P_p^S(T)=\sum_{m\geq 0} N_{p,m}(S)T^m \in \Z [[T]] \]

 The definition of the  \textit{geometric motivic Poincaré series},
 \[ \Pgeom^S(T)=\sum_{m\geq 0}\chi_c([j_m(H)])T^m \in \kmotq [[T]] \]
is inspired by that of the Serre-Oesterlé series.
However, there is no specialization of the series   $\Pgeom^S(T)$   into
 $ P_p^S(T)$ in general (see \cite{DL-Dwork}).

Denef and Loeser  studied the
 ``motivic nature'' of the series     $P_p^S(T)  $,  passing through
the Grothendieck ring $\kfield$  of ring formulas over $k$.
First, by Greenberg's theorem for every $m$ there exists a formula
$\psi_m$ over $k$ such that, for any field extension $k \subset K$, the $m$-jets over $k$
which can be lifted to arcs defined over $K$ correspond to the tuples satisfying
$\psi_m$ in $K$.  It follows that $\psi_m$ defines an element $[\psi_m] \in  \kfield$.
Then, Denef and Loeser defined a  ring
homomorphism $\chi_f:\kfield\rightarrow\kmot\otimes\Q$. This homomorphism
can be seen as a generalization of $\chi_c$,  since the image
by $\chi_f$ of the class of the ring formula defining a variety $V$ coincides
with the class $\chi_c([V])$ in $\kmot\otimes\Q$ .
The   {\em arithmetic motivic Poincaré series}  of $S$ is   defined as
\[\Par^S(T) =\sum_{m\geq 0}\chi_f([\psi_m])T^m  \in \kmotq [[T]].\]

Denef proved the rationality of the series $ P_p^S(T)$ using quantifier elimination results (see  \cite{Denef}).
Denef and Loeser proved the rationality of the series $\Pgeom^S (T)$ and $\Par^S(T)$
 by using quantifier eliminations theorems, various forms of motivic integration and
the existence of resolution of singularities (see \cite{DL-I, DL-JAMS}).

If  $V$ is a variety  defined over the integers and $p$ is a prime number,
the symbol $N_p (V)$ denotes the number of rational points of $V$ over the field of $p$ elements.
Denef and Loeser proved that the result of applying the operator $N_p$ to the motivic arithmetic series
$\Par^S (T)$ provides the Serre-Oesterlé series $P_p^S (T)$ for almost all  primes $p$.

If we fix the origin of the arcs in a fixed point  $0\in S$ we
obtain local versions of these series  $\Par^{(S, 0)} (T)$ and $\Pgeom^{(S,0)} (T)$,
which are also rational (see \cite{DL-I, DL-JAMS}).
The rationality proofs in \cite{DL-I, DL-JAMS} are qualitative in nature, in particular there is no conjecture on the
significance of the terms appearing in the denominator of the rational form of the series
$\Par^{(S,0)} (T)$ or in $\Pgeom^{(S,0)} (T)$.

The rational form of the series $\Par^{(S,0)}(T)$  is known explicitly for
a few classes of singularities. If $(S,0)$ is an analytically irreducible germ of plane curve,
the information provided by the series $\Par^{(S,0)}(T)$  is equivalent to the data of the
Puiseux pairs (see \cite{DL-JAMS}).
In \cite{Nicaise2} Nicaise proved the equality of the geometric and arithmetic motivic Poincaré series
in the case of varieties which admits a very special resolution of singularities,
in particular for normal toric surfaces (see also \cite{Nicaise1,LR}).
He gave a criterion for the equality $\Par^{(S,0)}(T) = \Pgeom^{(S,0)}(T)$  for various classes of singularities
and  also an example of normal toric threefold  $(S_0, 0) $ such that the series  $\Par^{(S_0, 0) }(T)$ and
$\Pgeom^{(S_0, 0)} (T)$ are different.
Some features of the motivic arithmetic
series are studied for quasi-ordinary singularities in \cite{Rond}.

In this paper we describe the arithmetic motivic Poincaré series of an affine toric variety
$Z^\Lambda = \mathrm{Spec} k[\Lambda]$, in terms of the semigroup $\Lambda$.
We assume that $\Lambda$ is a semigroup of finite type of a rank $d$ lattice $M$ (lattice of characters),  which
generates $M$ as a group, and such that the cone $\R_{\geq 0} \Lambda$ contains no lines.
In this situation there is a unique minimal system of generators $e_1, \dots, e_n$
of the semigroup $\Lambda$. The monomial ideal $(X^{e_i})_{i=1, \dots, n} \subset k [ \Lambda ]$ is maximal and
defines the distinguished point $0 \in Z^\Lambda$. In this paper we consider other monomial ideals as
the \textit{logarithmic jacobian ideals}   $\J_l$,
generated by monomials of the form $X^u$ for $u$ in the set
\[
\left\{ e_{i_1} + \dots + e_{i_l} \mid   e_{i_1} \y \dots \y e_{i_l} \ne 0 \right\}
\]
for  $l= 1, \dots, d$ (see  \cite{CoGP}),  and the \textit{ideals of sums of combinations} $\cc_j$, defined by
monomials $X^w$ with $w$ in the set
\[ \left\{ e_{i_1} + \dots + e_{i_j} \mid \{i_1,  \dots, i_j \}  \in    \binom{ \{1, \dots, n \} }{j}   \right\},
\]
where  $     \binom{ \{1, \dots, n \} }{j}$
denotes the set of combinations of $j$ elements of $\{1, \dots, n \}$, for $j=1, \dots, n$.

We study the motivic arithmetic series $\Par^{(Z^\Lambda, 0)} (T) $ by extending the approach
we used in
\cite{CoGP,CoGPqo} to describe the geometric motivic Poincaré series of toric and quasi-ordinary singularities.

By convenience we explain the methods and results first when the variety $Z^\Lambda$ is normal.
The set $j_m (H (Z^\Lambda)_0 )$ of $m$-jets of arcs through $(Z^\Lambda, 0)$
is constructible; it is  a finite disjoint union
of locally closed subsets of the form $j_m (H^*_\nu)$ (see \cite{CoGP}). Here $H^* _\nu$ denotes the set of arcs through    $(Z^\Lambda, 0)$
which have generic point in the torus and a given order $\nu \in M^*$. The set $H^*_\nu$ is an orbit of the
natural action of the arc space of the torus on the arc space of the toric variety $Z^\Lambda$
(see \cite{Ishii-algebra, Ishii-crelle}).

We describe the class, denoted by $\chi_f ([j_m (H^* _\nu)]_f ) $,  of the formula
defining the locally closed subset $j_m (H^*_\nu)$ in terms of
the Newton polyhedra of the logarithmic jacobian ideals
and  the degree of certain  Galois cover.  This Galois cover reflects the relation
between the  initial coefficients of the arcs in $ H^*_\nu $ and
the initial coefficients of the $m$-jets in   $j_m (H^*_\nu)$, see Section  \ref{coeffs}.

A key point in the description of the rational form of the series  $\Par^{(Z^\Lambda, 0)} (T) $
is that
using the ideals $\cc_j$ we can refine a finite partition of the set of possible pairs $\{ (\nu, m) \}$, which
was defined in \cite{CoGP} to describe the sum of $\Pgeom^ {(Z^\Lambda, 0)} (T) $. If   $(\nu, m)$ and
$(\nu', m')$ belong to the same subset of this refinement then the degrees of the Galois covers associated to
$j_m  (H^*_\nu)$ and $j_{m'}(H^*_{\nu'} )$ coincide (see Sections \ref{fan} and \ref{Combinatoria}).
Using these partitions we decompose the series $\Par^{(Z^\Lambda, 0)} (T) $
as a sum of a finite number of contributions.
The main result is a formula for the rational form of $\Par^{(Z^\Lambda, 0)} (T) $
(see Theorem \ref{ThNormal} and Corollary \ref{P-geom}).
The proofs pass by the results on the generating function
of the projection of the set of integral points in the interior of a rational polyhedral cone (see  \cite{CoGP}).
The denominator of $\Par^{(Z^\Lambda, 0)} (T) $ is a finite product of terms of the form $1 - \L^a T^b$ with
 $a \geq 0$  and $b >0$, which are determined explicitly in terms of the ideals of sums of combinations $\cc_j$.
The integers $a$ and $b$ can be described in terms of the orders of vanishing of the ideals $\cc_j$ and $\J_l$
at the codimension one orbits of various toric modifications given by the Newton polyhedra of
the ideals $\cc_j$ (see Remark \ref{inter}).
In the normal toric case we obtain a formula for $\Par^{Z^\Lambda} (T)$ in terms of
 arithmetic motivic series at the distinguished points of the orbits.

In the non-normal case, we obtain in a similar way a formula for the rational form of
 $\Par^{(Z^\Lambda, 0)} (T) $ and the factors of
its denominator. The main difference is that we have to consider contributions
of jets of arcs with generic point in the various orbits of $Z^\Gamma$.
                 We deduce a formula for the difference $\Pgeom^{(Z^\Lambda, 0)} (T) - \Par^{(Z^\Lambda, 0)} (T)$
and we give a criterion for the equality of these two series which generalizes the one given by Nicaise in
\cite{Nicaise1} (see Proposition \ref{compara2} and Theorem \ref{Pg=Pa}).

The paper is organized as follows.
In Sections \ref{add-inv} and \ref{ArcsJetsMotseries} we introduce
the Grothendieck rings, the arc and jet spaces and the
motivic Poincaré series.
The notations on toric varieties, their monomial ideals and their arcs are introduced in   Sections \ref{sec-tor}
and \ref{ArcsAndJets}.
The computation of the class  $\chi_f ([j_m (H^* _\nu)]_f ) $ is given in  Section \ref{coeffs}.
Sections \ref{fan} and  \ref{Combinatoria}
deal with the partitions associated to sequences of monomial ideals.
The main results are stated and proved in Sections  \ref{ArithToric},    \ref{RatGenSeries}  and \ref{main}.
In the case of normal toric varieties some features of the computation can be simplified (see Section \ref{NormalCase}).
We discuss some examples in Section \ref{ejemplos}.

\section{Grothendieck rings of varieties and of ring formulas}
\label{add-inv}
The \textit{Grothendieck ring $\kvar$ of $k$-varieties} is the
free abelian group of isomorphism classes $[X]$ of $k$-varieties
$X$ modulo the relations $[X]=[X']+[X\setminus X']$ if $X'$ is
closed in $X$, and  where the product is defined by
$[X][X']=[X\times  X']$. We denote by $\L:= [\A_k^1]$ the class of
the affine line.
If $C$ is a constructible subset of some variety $X$,
i.e. a disjoint union of finitely many locally closed subvarieties
$A_i$ of $X$, then  $[C] \in \kvar$ is well defined as $[C]
:=\sum_i [A_i ]$ independently of the representation. Bittner
proved, using the weak factorization theorem, that the ring
$\kvar$ is generated by classes of smooth projective
$k$-varieties, modulo relations of the form $[W] - [E] = [X] -
[Y]$, where $Y \subset X$ is a closed subvariety, and $W$ is the
blowing up of $X$ along $Y$ with exceptional divisor $E$ (see  \cite{ Bittner}).

There exists a unique ring homomorphism:
\begin{equation}
 \label{chi}
\chi_c : \kvar \to \kchow,
\end{equation}
which maps the class of a smooth projective variety over $k$ to
its \textit{Chow motive}, where $\kchow$ denotes the
\textit{Grothendieck ring of the category of Chow motives} over
$k$ (with coefficients in $\Q$). This result, which is due to
Guillet and Soulé  \cite{GS} and Guillén and Navarro Aznar
\cite{GN}, can be seen also in terms of Bittner's  result.   We
refer to \cite{GS, GN, Bittner} for details and to \cite{Scholl}
for an introduction to the notion of motives.  We denote by
$\kmot$ the image of $\kvar$ in $\kchow$ under this homomorphism.
Notice that the image of $\L$ in $\kmot$, which we denote with the
same symbol, is not a zero divisor in $\kmot$ since it is a unit
in $\kchow$.

A \textit{ring formula} $\psi$ over $k$ is a first order formula
in the language of $k$-algebras and free variables
$x_1,\ldots,x_n$, that is,  the formula $\psi$ is built from
boolean combinations ("and", "or", "not") of polynomial equations
over $k$ and existential and universal quantifiers. The
\textit{Grothendieck ring $\kfield$ of ring formulas over $k$}, is
generated by symbols $[\psi]$, where $\psi$ is a ring formula over
$k$, subject to the relations
$[\psi_1\vee\psi_2]=[\psi_1]+[\psi_2]- [\psi_1\wedge\psi_2]$ if
$\psi_1$ and $\psi_2$ have the same free variables, and
$[\psi_1]=[\psi_2]$ if there exists a ring formula $\Psi$ over $k$
such that, when interpreted in any field $K\supseteq k$ provides
the graph of a bijection between the tuples of elements of $K$
satisfying $\psi_1$ and those satisfying $\psi_2$. The ring
multiplication is induced by the conjunction of formulas in
disjoint sets of variables (see \cite{DL-JAMS}). Denef and Loeser
defined a ring homomorphism
\begin{equation}   \label{chi-formula}
        \chi_f: \kfield \rightarrow \kmot \otimes\Q.
\end{equation}
They proved that this homomorphism is characterized by two
conditions. The first one is that for any ring formula $\psi$
which is a conjunction of polynomial equations over $k$, the
element $\chi_f([\psi])$ is equal to the class $\chi_c ([V])$  in
$\kmot \otimes\Q$ of the variety  $V$ defined by $\psi$. The
second condition, which is more technical, expresses that certain
relations should hold in terms of unramified Galois coverings over
$k$. We refer to \cite{DL-JAMS,DL-Dwork} for the precise
statement. In the simplest case it implies the following:

\begin{Exam}   (see \cite{DL-Dwork} Example 6.4.3)
       If $n \geq 1$ is a fixed integer,  $k$ is a field containing all $n$-th roots of unity and
     $\psi$ is the ring formula $\psi:(\exists y)(x=y^n\mbox{ and }x\neq 0)$
then we have that $\chi_f([\psi])=\frac{1}{n}(\L-1)$.
\end{Exam}

We deduce from this example the following Lemma:
\begin{Lem}
Let $\psi$ be the ring formula whose interpretation in any field
$K\supseteq k$ provides the set of $K$-rational points of $T$
which lift to $K$-rational points of $T'$ by a Galois covering
$T'\rightarrow T$ of degree $n$ of $d$-dimensional algebraic
$k$-tori. If the field $k$ contains all the $n$-th roots of unity
then we have that $\chi_f([\psi])=\frac{1}{n}(\L-1)^d$.
\label{Lema0}
\end{Lem}

{\em Proof.} The morphism $T'\rightarrow T$ induces a finite index inclusion of the corresponding character group
$M\subseteq M'$, and hence a map of $k-$algebras $k[M]\hookrightarrow k[M']$. By the classification theorem of finitely
generated abelian  groups applied to $M'/M$ there exists a basis $\{v_1,\ldots,v_d\}$ of $M'$ and unique integers
$b_1|b_2|\cdots|b_d$, where $|$ denotes division, such that $\{b_1v_1,\ldots,b_dv_d\}$ is a basis of $M$ and
$n=b_1\cdots b_d$. It follows that the map of coordinate rings $K[M]\hookrightarrow K[M']$ express in coordinates as
$K[z_1^{±b_1},\ldots,z_d^{±b_d}]\hookrightarrow K[z_1^{±1},\ldots,z_d^{±1}]$.
 We deduce that the ring formula $\psi$ is the conjunction of formulas
$\psi_i:(\exists y_i)(x_i=y_i^{b_i}\mbox{ and }x_i\neq 0)$, for
$i=1,\ldots,d$ where the variables $x_1,\ldots,x_d$ are
independent. Then we get that $ \chi_f([\psi])=\frac{1}{b_1\cdots
b_d}(\L-1)^d$. \hfill $\Box$

\begin{remark}
Denef and Loeser defined the map $\chi_f$ by factoring it through the  Grothendieck ring $\kpff$
of ring formulas for the category of \textit{pseudo-finite fields} containing $k$.
See \cite{DL-JAMS, DL-Dwork, DL-ICM}.
\label{remAss2}
\end{remark}

\section{Arcs, jets spaces and motivic Poincaré series}
\label{ArcsJetsMotseries}

We start this Section  by recalling the definition of the space of
arcs of a variety $S$. We assume for simplicity that $S$ is an
affine  irreducible and reduced algebraic variety  defined over a
field $k$ of characteristic zero.

For any integer $m \geq 0$ the functor from the
category of $k$-algebras to the category of sets, sending a
$k$-algebra $R$ to the set of $R[t]/(t^{m+1})$-rational points of
$S$ is representable by a $k$-scheme $H_m(S)$ of finite type over
$k$, called the $m$-jet scheme of $S$. The natural maps induced
by truncation $j_m^{m+1}: H_{m+1} (S) \rightarrow H_m (S)$ are
affine and hence the projective limit $H(S):= \limproj H_m (S)$ is
a $k$-scheme, not necessarily of finite type, called the \textit{arc
space} of $S$.

In what follows we consider the schemes $H_m(S)$ and $H(S)$ with
their reduced structure.     We have natural morphisms
$j_m: H(S) \rightarrow H_m (S)$. By an {\em arc} we mean a
$k$-rational point of $H(S)$, i.e., a morphism $ \mbox{Spec }
k[[t]] \rightarrow S$. By an $m$-jet we mean a $k$-rational
point of $H_m(S)$, i.e., a morphism $ \mbox{Spec } k[t]/(t^{m+1})
\rightarrow S$. The origin of the arc (resp. of the $m$-jet) is
the image of the closed point $0$ of $ \mbox{Spec } k[[t]]$
(resp. of $ \mbox{Spec } k[t]/(t^{m+1}) $).

            If $Z \subset S$ is a closed
subvariety then $H(S)_Z:= j_0 ^{-1} (Z) $ (resp. $H_{m}(S) _Z:=
(j^m_0 )^{-1} (Z)$) denotes the subscheme of $H(S)$ (resp. of
$H_{m}(S)$) formed by arcs (resp. $m$-jets) in $S$ with origin in
$Z$.

 By a Theorem of Greenberg
\cite{Gr}, for any integer $m \geq 0$,  $j_m (H(S))$ is a
constructible subset of the $k$-variety  $H_m (S)$.     We can
then consider the class   $[j_m (H(S))] \in \kvar$. Greenberg's
result implies also  that there is a ring formula $\psi_m$ over
$k$, such that for any field $K$ containing $k$, the $k$-rational
points of $H_m (S)$ which can be lifted to $K$-rational points of
$H(S)$  correspond to the tuples satisfying $\psi_m$ in $K$. If
$\psi_m'$  is another ring formula over $k$ with the same property
then $[\psi_m] = [\psi_m']$ in $\kfield$.
 The same applies for $j_m (H(S)_Z)$ if $Z\subset S$ is a closed
  subvariety.

\begin{Not}
      We denote the class $[\psi_m]$ by $[j_m (H(S))]_f$  to avoid confusion with
 the class  $[j_m (H(S))] \in \kvar$.
\end{Not}

      The following Poincaré series were introduced by Denef and Loeser in the papers (\cite{DL-I,DL-JAMS}).
 \begin{definition}  $\,$      \label{par-geom}
\begin{enumerate}
 \item        The \textit{geometric motivic Poincaré series} of $(S,Z)$
is
\[
 \Pgeom^{(S, Z)} (T) := \sum_{m \geq 0} \chi_c ([ j_m (H(S)_Z)
]) T^m \in \kmotq [[T]].
\]

\item  The \textit{arithmetic motivic Poincaré series}  of $(S,Z)$ is
 \[\Par^{(S, Z)} (T) := \sum_{m \geq 0} \chi_f ([ j_m (H(S)_Z)]_f ) T^m \in \kmotq [[T]]. \]
\end{enumerate}
\end{definition}

\begin{remark}
We have slightly modified  the original definition of the
 geometric motivic Poincaré series, as   $\sum_{m \geq 0} [ j_m (H(S)_Z) ]  T^m  \in \kvar [[T]]$ (see \cite{DL-I}),
 in order to have the geometric and arithmetic setting
in the same ring. This does not affect the rationality results below.
\end{remark}

Denef and Loeser proved that these series have a rational form:

\begin{The} (see \cite{DL-I}   Theorem 1.1 and  \cite{DL-JAMS} Theorem 9.2.1)  The series $\Pgeom^{(S,Z)} (T)$ and
$\Par^{(S, Z)} (T)$   belong to
the subring of $\kmotq [[T]]$ generated by $\kmotq[T]$ and the series
$(1-\L^aT^b)^{-1}$, with $a\in\Z$ and  $b>0$.
\label{RacDL}   \label{RatDLArit}
\end{The}

The arithmetic motivic  Poincaré series has interesting properties of specialization to classical arithmetic series.
Let $p$ be a prime number. The operators   $N_p$ and $N_{p, m}$ are  applied to a variety $V$ defined over the integers by
$
 N_p (V) := \# V( \Z / p\Z )$   and $N_{p, m} (V) := \# \{ \p_m (V(\Z_p))  \}$
where $\Z_p$ denotes the $p$-adic integers, $ \p_m (V(\Z_p))
\subset V(\Z / p^{m+1} \Z) $ is the image of $V(\Z_p))$ by the
natural projection induced by $\Z_p \to \Z/ p^{m+1}\Z$, and $\#$
denotes the cardinality. Suppose that the variety $S$ is defined
over the integers. The  \textit{Serre-Oesterlé} series $P_p ^S
(T)  := \sum_{m \geq 0} N_{p, m} T^m \in \Z[[T]]$ of $S$ at the
prime $p$
 is a rational function of $T$ (see \cite{Denef}).
Denef and Loeser proved that for $p \gg 0$ the series $P_p^S (T)$
is obtained from $\Par^S (T)$ by applying to each coefficient the
operator $N_p$ (see \cite{DL-JAMS, DL-ICM, DL-Dwork}).

\begin{remark}
These results hold in a more general setting,  in particular when
$S$ is not affine  as assumed here (see \cite{DL-I,
DL-JAMS}).    The proof of the rationality of  $P_p ^S (T)$
involves the use of quantifier elimination results and $p$-adic
integration (see \cite{Denef}).  The proof of the  rationality of
$\Par^S (T)$ requires also  quantifier elimination
results and arithmetic motivic integration (see
\cite{DL-JAMS, DL-Dwork, DL-ICM}).

\end{remark}

\section{Affine toric varieties and monomial ideals} \label{sec-tor}

In this Section we introduce  the basic notions and notations
from toric geometry(see \cite{Ewald,Oda,Fu,GKZ}
  for
proofs).

If $N \cong \Z^{d}$ is a lattice we denote by $N_\R:=N\otimes\R$
(resp. $N_\Q:=N\otimes\Q$) the vector space spanned by $N$ over
the field $\R$ (resp. over $\Q$). In what follows a {\em cone} in
$N_\R$ mean a {\em rational convex polyhedral cone}: the set of
non negative linear combinations of vectors $a_1 \dots , a_r \in
N$. The cone $\t$ is {\em strictly convex} if it contains no line
through the origin,  in that case we denote by $0$ the
$0$-dimensional face of $\t$;
the cone $\t$ is {\em simplicial} if the primitive vectors of the
$1$-dimensional faces are linearly independent over $\R$. We
denote by $\stackrel{\circ}{\t}$ or by $\mbox{{\rm int}} (\t)$ the
relative interior of the cone $\t$.

 We denote by $M$ the dual
lattice. The {\em dual} cone  $\t^\vee \subset M_\R$ (resp. {\em
orthogonal} cone $\t^\bot$) of $\t$ is the set $ \{ w  \in M_\R\ |\,
\langle w, u \rangle \geq 0,$  (resp. $ \langle w, u \rangle = 0$)
$ \; \forall u \in \t \}$.

A {\em fan} $\Sigma$ is a family of strictly convex
  cones  in $N_\R$
such that any face of such a cone is in the family and the
intersection of any two of them is a face of each. The relation
$\theta \leq \t$  (resp. $\theta < \t$) denotes that $\theta$ is a
face of $\t$ (resp. $\theta \ne \t$ is a face of $\t$). The {\em
support} (resp. the $k$-{\em skeleton}) of the fan $\Sigma$ is the
set $|\Sigma | := \bigcup_{\t \in \Sigma} \t \subset N_\R$ (resp.
$\Sigma^{(k)} = \{ \t \in \Sigma \mid \dim \t = k \}$).
We say that a fan $\Sigma'$ is a {\em subdivision\index{fan
subdivision}} of the fan $\Sigma$ if both fans have the same
support and if every cone of $\Sigma'$ is contained in a cone of
$\Sigma$. If $\Sigma_i$ for $i=1,\dots, n$, are fans with the
same support their intersection
$\cap_{i=1}^n\Sigma_i:=\{\cap_{i=1}^n\t_i\ |\,
\t_i\in\Sigma_i\}$ is also a fan.

\begin{Not} \label{Lambda}
In this paper $ \Lambda  $ is a sub-semigroup of finite type of a
lattice $M$,  which generates $M$ as a group and such that the
cone $\s^\vee = \R_{\geq 0} \Lambda$ is strictly convex and of
dimension $d$. We denote by $N$ the dual lattice of $M$ and by $\s
\subset N_{\R}$ the dual cone of $\s^\vee$. We denote by
$Z^\Lambda$ the {\em affine toric variety} $Z^{\Lambda} =
\makebox{Spec} \, {k}[ \Lambda ]$, where $ {k} [ \Lambda] = \{
\sum_{\rm finite}  a_\l{X}^\l\mid a_{\l} \in {k}\} $ denotes the
semigroup algebra of the semigroup $\Lambda$ with coefficients in
the field ${k} $. The semigroup $\Lambda$ has a unique minimal set
of generators $e_1,\dots,e_n$ (see the proof of  Chapter V, Lemma
3.5, page 155 \cite{Ewald}). We have an {\em  embedding} of
$Z^\Lambda\subset \A_{k} ^n$ given by, $x_i := X^{e_i}$ for $i=1,
\dots,  n$.
\end{Not}

If $\Lambda = \s^\vee \cap
M$   then the variety    $Z^{\Lambda}$,
which we denote also by         $Z_{\s, N}$ or by
 $Z_\s$ when the lattice is clear from the context, is normal.
If $\Lambda \ne \s^\vee \cap
M$       the
inclusion of semigroups $\Lambda \rightarrow \bar{\Lambda} $
defines a toric modification $Z^{\bar{\Lambda}} \rightarrow
Z^{{\Lambda}}$,
  which  is the {\em normalization map}.

The torus $ T_N:= Z^{M}$ is an open dense subset of $Z^\Lambda$,
which acts on $Z^\Lambda$ and the action extends the action of the
torus on itself by multiplication. The {\em origin} $0$ of the
affine toric variety $Z^{\Lambda}$ is the $0$-dimensional orbit,
defined by the maximal ideal $(X^{\l})_{0 \ne \l\in
\Lambda }$ of ${k}[\Lambda]$.
 There is a one to one inclusion
 reversing correspondence between the faces  of $\s$ and the orbit
 closures  of the torus action on $Z^\Lambda$.
 If $\theta \leq \s$,
 we denote by  $\O^\Lambda_\theta$ the orbit corresponding to the
 face $\theta$ of $\s$.
The orbit closures are of the form  $Z^{\Lambda \cap \theta^\bot }$
for $\theta \leq \s$.

The {\em Newton polyhedron} of a monomial ideal corresponding to a
non empty set of lattice vectors
 ${\mathcal I} \subset   \Lambda $ is defined as the convex hull
of the Minkowski sum  of sets ${\mathcal I} + \s^\vee$. We denote
this polyhedron by ${\mathcal N} ({\mathcal I})$.    Notice
that the vertices of ${\mathcal N} (\mathcal I)$ are elements of
$\mathcal I$.
We denote by $\mbox{\rm ord}_{\mathcal {\mathcal I}}$ the {\em support
function} of the polyhedron ${\mathcal N} ( {\mathcal I} )$, which
is defined by $\mbox{\rm ord}_{\mathcal I} : \s \rightarrow \R$, $
\nu \mapsto \inf_{\omega \in {\mathcal N} ( {\mathcal I} )}
\langle \nu, \omega \rangle$.
The face  of  the polyhedron ${\mathcal N} ({\mathcal I})$
determined by $\nu \in \s$ is  the set $ {\mathcal F}_\nu := \{
\omega \in {\mathcal N} (\mathcal I) \mid \langle \nu, \omega
  \rangle   = \mbox{\rm ord}_{\mathcal I} (\nu ) \}$.
All faces of ${\mathcal N} ({\mathcal I})$ are of this form, the
compact faces are defined by vectors $\nu  \in
\stackrel{\circ}{\sigma}$. The set   $\Sigma ({\mathcal
I})$    consisting  of the cones $ \s(
{\mathcal F} ) := \{ \nu \in \s \; \mid \langle \nu , \omega
\rangle   = \mbox{\rm ord}_{\mathcal I} ( \nu ),
 \; \forall \omega \in {\mathcal F}\}$,
for ${\mathcal F}$ running through the faces of ${\mathcal N}
({\mathcal I})$, is a fan supported on $\s$.
Notice that if $\t\in  \Sigma ( {\mathcal I} )$ and if
$\nu, \nu' \in \stackrel{\circ}{\t}$ then ${\mathcal
F}_{\nu} = { \mathcal F}_{\nu'}$. We denote this face of
${\mathcal N} (\mathcal I)$ also by  ${ \mathcal F}_{\t}$.

The
 affine varieties $Z_\t$ corresponding to cones $\t$ in a fan $\Sigma$
 glue up to define a {\em toric variety} $
 Z_\Sigma$.
A fan  $\Sigma$ subdividing the cone $\sigma$ defines a {\em toric
 modification} $ \pi_{\Sigma} : Z_{\Sigma}
 \rightarrow   Z_\s$.

If   ${\mathcal I} \subset   \Lambda $ defines a monomial ideal the composite
  $   Z_{ \Sigma ({\mathcal I})}  \stackrel{ \p_{\Sigma({\mathcal I})}} {\longrightarrow} Z_\s \longrightarrow  Z^\Lambda$
is equal to the {\em normalized blowing up} of $Z^\Lambda$ centered at ${\mathcal I}$ (see \cite{LR} for
  instance).

\begin{definition}
For $1 \leq j \leq n$ the  {\em $j$-th ideal of sums of
combinations} of $Z^\Lambda$ is the monomial ideal $\cc_j$ of ${k}
[\Lambda]$ generated by $X^\a$ where $\a$ runs through:
\begin{equation}       \label{c-j}
\left\{ e_{i_1} + \dots + e_{i_j} \mid \{i_1,  \dots, i_j \}  \in    \binom{ \{1, \dots, n \} }{j}   \right\},
\end{equation}
where  $     \binom{ \{1, \dots, n \} }{j}$
denotes the set of combinations of $j$ elements of $\{1, \dots, n \}$, for $j=1, \dots, n$.
 We denote by $\Theta_j$  (resp.
by $\ord_{\cc_j}$) the dual subdivision of $\s$ (resp. the support
function) of the polyhedron $\Newton (\cc_j) $. The maps
\[
\begin{array}{lcccccccl}
\varphi_1    &   :=  &   \mbox{\rm ord}_{\cc_1  } &  \mbox{ and }
& \varphi_j  & :=  &   \mbox{\rm ord}_{\cc_j }  - \mbox{\rm
   ord}_{\cc_{j-1}  } &  \mbox{ for }  j=2,\dots,n,
\end{array}
\]
are piece-wise linear functions defined on the cone $\s$. If $\nu
\in \s$ we denote by convenience $\varphi_0 (\nu) := 0$ and
$\varphi_{n+1} (\nu): = + \infty$.

\label{defIdeals}
\end{definition}

\begin{definition} For $1 \leq l \leq d$ the  {\em
$l$-th logarithmic jacobian ideal} of $Z^\Lambda$ is the monomial
ideal $\J_l$ of ${k} [\Lambda]$ generated by $X^\a$ where $\a$
runs through:
\begin{equation}       \label{j-k}
\{ e_{i_1} + \cdots + e_{i_l}  \;  \mid \; e_{i_1} \y \cdots \y
e_{i_l}  \ne 0,
  \mbox{ for }
1 \leq i_1 , \dots, i_l \leq n  \}.
\end{equation}
We denote by $\Sigma_l $  (resp. by $\ord_{\mathcal{J}_l }$) the
dual subdivision of $\s$ (resp. the support function) of the
polyhedron  $\Newton ({\mathcal J}_l) $. The maps
\[
\begin{array}{lcccccccl}
\f_1    &   :=  &   \mbox{\rm ord}_{\mathcal J_1  } &  \mbox{ and
} & \f_l  & :=  &   \mbox{\rm ord}_{\mathcal J_l  }  - \mbox{\rm
   ord}_{\mathcal J_{l-1}  } &  \mbox{ for }  l=2,\dots,d,
 \end{array}
\]
are piece-wise linear functions defined on the cone $\s$. If $\nu
\in \s$ we denote by convenience $\f_0 (\nu) := 0$ and $\f_{d+1}
(\nu): = + \infty$.
\end{definition}

We  use the notation $\J_l$ (resp. $\cc_j$) also for the set
(\ref{j-k}) (resp. (\ref{c-j})).

\begin{Lem} \label{ine} If $\nu \in
\stackrel{\circ}{\s}$ and if $(p_1, \dots, p_n)$ is  a permutation
of $(1, \dots, n)$ such that
\[
\langle \nu, e_{p_1} \rangle \leq \cdots \leq \langle \nu, e_{p_n}
\rangle.
\]
then  $\ord_{\cc_j} (\nu) = \langle \nu, \sum_{r=1}^j e_{p_r}
\rangle$ and $\varphi_j (\nu) = \langle \nu, e_{p_j}\rangle$ for
$1 \leq j \leq n$. Moreover, the following holds
\[
0 = \varphi_0  (\nu) \leq \varphi_1 (\nu) \leq \dots\leq
\varphi_{n} (\nu) \quad \mbox{ and } \quad  0 =   \f_0  (\nu) \leq
\f_1 (\nu) \leq \dots\leq \f_{d} (\nu).
\]
\end{Lem}
{\em Proof.} The first assertion follows by induction on $j\in\{1,\ldots,n\}$.

 See Lemma 5.3 \cite{CoGP} for the second sequence of inequalities.
\hfill $\ {\Box}$

\begin{Pro}
The Newton polyhedra of the ideals $\mathcal C_j$, $j=1,\ldots,n$,
determine and are determined by the minimal system of generators of the semigroup $\Lambda$.
\label{ObsSemigrupo}
\end{Pro}
{\em Proof.}
The Newton polyhedron $\mathcal N(\mathcal C_j)$ determines and it is determined by
its support function  $\ord_{\mathcal C_j}$, for
$j = 1, \dots, n$.
By Lemma \ref{ine} and the definitions
if $\theta$ is a $d$ dimensional cone of the fan $\cap_{r=1} ^n \Theta_r$
there exists a permutation $i_1, \dots, i_n$ of $1, \dots, n$ such that
$
\varphi_j (\nu) = \langle \nu, e_{i_j}  \rangle
$ for $j=1,  \dots, n$ and all
 $\nu \in \stackrel{\circ}{\theta}$.
Thus, the functions  $\varphi_j$, $j=1, \dots, n$, , or equivalently,
$\ord_{\mathcal C_j}$,
$j = 1, \dots, n$,
determine the vectors $e_1, \dots, e_n$.
\hfill $\Box$

\section{Arcs and jets on a toric singularity}

\label{ArcsAndJets}

Let $\Lambda$ be a semigroup as in Notation \ref{Lambda}.
If $R$ is a $k$-algebra, a $R$-rational point
of $Z^{\Lambda}$ is a homomorphism of semigroups $(\Lambda, +)
\rightarrow (R,\cdot)$, where $(R,\cdot)$ denotes the semigroup
$R$ for the multiplication. In particular, the closed points are
obtained for $R = k$.
An arc
$h$ on the affine toric variety $Z^{\Lambda}$ is given by a semigroup
homomorphism $(\Lambda, +) \rightarrow (k[[t]], \cdot)$. An
arc in the torus $ T_N$ is defined by a semigroup homomorphisms $\Lambda
\rightarrow k[[t]]^*$, where $k[[t]]^*$ denotes the group of
units of the ring $k[[t]]$.

\begin{Not}    We denote  the set of arcs $H(Z^\Lambda)_0$
of $Z^\Lambda$ with
origin at the distinguished point $0$ of $Z^\Lambda$  simply by $ H_\Lambda$, and by
    $H_\Lambda^*$ the set consisting of those arcs of
$H_\Lambda$ with generic point in the torus $ T_N$.
\label{NotacionArcs}
\end{Not}
Notice that $h \in H_\Lambda^*$ if and only if for all $u \in
\Lambda$ the formal power series ${X}^u \circ h \in
k [[t]]$ is non-zero.
Any arc  $h  \in H^*_\Lambda$ defines two group
 homomorphisms $\nu_h : M \rightarrow \Z  \mbox{ and } \omega_h: M
\rightarrow k[[t]]^* \mbox{ by: } {X}^m \circ {h} = t^{\nu_h (m)}
\omega_h (m)$.  If $m \in \Lambda$ then  $\nu_h (m)
>0$ hence $\nu_h$ belongs to $\stackrel{\circ}{\s} \cap N$.
Notice that $\omega_h$ defines an arc in the torus, i.e.,
$\omega_h \in H ({ T_N}) $.

\begin{remark}
The space of arcs in the torus acts on the arc space of a toric
variety (see \cite{Ishii-algebra, Ishii-crelle}).
\end{remark}

\begin{Lem} \label{encode} {\rm (see Theorem 4.1 of \cite{Ishii-algebra}, Lemma 5.6 of
 \cite{Ishii-crelle}, and Proposition 3.3 \cite{LR}).}
The map $ \stackrel{\circ}{\s} \cap N \times  H({ T_N}) \rightarrow
H^*_\Lambda$  which applies a pair $(\nu,\omega)$ to the arc $ h$
defined by
$
{X}^u \circ {h} = t^{\langle \nu, u \rangle} \omega (u), \mbox{
for } u \in \Lambda ,
$
is a one to one correspondence. The sets $H^*_{\Lambda, \nu} := \{
h \in H^*_\Lambda\ |\ \nu_h = \nu \}$   for $\nu \in
\stackrel{\circ}{\s} \cap N$ are orbits for the action of $H_{
T_N} $ on $H^*_\Lambda$ and we have that $H_\Lambda^* =
\bigsqcup_{\nu \in  \stackrel{\circ}{\s} \cap N}  H^*_{\Lambda,
\nu} $.
\end{Lem}

\begin{remark}
We often denote the set $H^*_\Lambda$ (resp. the orbit
$H^*_{\Lambda, \nu} $) by  $H^*$ (resp. by $H^*_\nu$) if $\Lambda$
is clear from the context.
\end{remark}

 An arc $h \in H_{\Lambda}$
has its generic point $\eta$ contained in exactly one orbit of the
torus action on $Z^\Lambda$.  If    $h(\eta) \in
\O_\theta^\Lambda$, for some $\theta \leq \s$, then $h$ factors
through the orbit closure $Z^{\Lambda \cap \theta^\bot}$ and  $h
\in H_{\Lambda \cap \theta^\bot}^*$, i.e., $h$ is an arc through
$(Z^{\Lambda \cap \theta^\bot}, 0)$ with generic point in the
torus  $\O_\theta^\Lambda$.
     We can apply Lemma \ref{encode} to describe the set
$H_{\Lambda \cap \theta^\bot}^*$, just replacing the semigroup
$\Lambda$ by $\Lambda \cap \theta^\bot$ (see \cite{CoGP}).
 In
particular, if $\theta = 0$ then $h \in H_\Lambda^*$;  if $\theta
=\s$ then $\Lambda \cap \theta^\bot =0$ and $h$ is the constant
arc at the distinguished point $0 \in Z^\Lambda$.
We have a
partition $     H_\Lambda = \bigsqcup_{\theta \leq \s}  H_{\Lambda \cap
    \theta^\bot}^*$.

\section{The image of the class of the formula defining $j_s (H^*_\nu)$}
\label{coeffs}

\begin{definition} We associate to  $(\nu, s) \in (\stackrel{\circ}{\s} \cap
N) \times  \Z_{>0}$ the sets
\[
M_\nu^s :=  \span_\Z \{ e_i \mid \langle \nu, e_i\rangle \leq s, i
= 1, \dots, n \}, \quad \mbox{ and }  \quad  \ell_\nu^s :=
\span_\Q \{ e_i \mid \langle \nu, e_i\rangle \leq s, i = 1, \dots,
n \}.
\]
We denote by $l (\nu, s)$ the dimension of the $\Q$-vector space
$\ell_\nu^s$. The integer $l (\nu, s)$ is also the rank of the
lattice $M_\nu^s$. We denote by $q(\nu,s)$ the index of the
lattice extension $M_\nu^s \subset \ell_\nu^s \cap M$.
\label{defLattices}
\end{definition}

\begin{Pro}   If  $(\nu,s)\in \stackrel{\circ}{\s} \times  \Z_{> 0}$, $l(\nu, s) > 0$, and if the field $k$ contains
all the $q(\nu,s)$-th roots of unity then we have
\[ \chi_f ([j_s(H_\nu^*)]_{f}) = \frac{1}{q(\nu,s)}  (\L-1)^{l(\nu,s)} \times  \L^{sl (\nu,s)
 -\ord_{\J_{l(\nu,s)}}(\nu)}. \]
 If $l(\nu, s) = 0$ then we have  $\chi_f ([j_s(H_\nu^*)]_{f}) =
 1$.
\label{ProArithKey}
\end{Pro}

{\em Proof.}  If $h \in H^*_\nu $   the equality  $ \ord_t
 ( X^{e_i} \circ h ) = \langle \nu, e_i \rangle $ holds for $1 \leq i \leq n$.
By Definition    \ref{defLattices} those vectors $e_i$ such that
$j_s( X^{e_i} \circ h) \ne 0$ span the     $\Q$-vector space
$\ell_\nu^s$ since $ \langle \nu, e_i \rangle \leq s$.
 If $l(\nu,s)=0$ this vector space is empty, the jet space $j_s (
H^*_{\nu})$ consists of the constant $0$-jet and the conclusion
follows easily from the definitions.

Suppose then that $l:= l(\nu,s ) >0$. If $h\in H_\nu^*$ then it is
given by $n$ series of the form
\[
X^{e_i}\circ h=t^{\langle\nu,e_i\rangle}c(e_i)\big(1+ \sum_{m \geq
1} u_m (e_i)t^m \big), \ i=1,\ldots,n. \] We have that the $s$-jet
$j_s (X^{e_i} \circ h)$ is different from zero if and only if
$\langle \nu, e_i \rangle \leq s$.

By Lemma 5.7 of \cite{CoGP} there exist integers $ 1 \leq
{k_1},\dots, {k_l} \leq n $ such that $\phi_i (\nu) = \langle
\nu, e_{k_i} \rangle \leq s$, for $i=1, \dots, l$, $\ell_\nu^s  =
\span_\Q \{e_{k_1}, \dots, e_{k_l} \}$ and $\ord_{\J_l} (\nu) =
\sum_{i=1}^l \langle \nu, e_{k_i} \rangle$.

 By Section 6 of
\cite{CoGP} if $h$ is the universal family of arcs   parametrizing $H^*_\nu$,
the terms $\{ u_m (e_{k_i}) \mid i=1, \dots, l, m \geq 1 \}$ are
algebraically independent over $\Q$ and the terms $\{ c(e_i)^{±
1} \mid i=1, \dots, n \}$ generate a $k$-algebra isomorphic to $k
[M]$ by the isomorphism which maps $c(e_i) \mapsto X^{e_i}$.

By the proof of Theorem 7.1 \cite{CoGP} a formula defining $j_s(
H^*_\nu)$ is the conjunction of two formulas $\psi_1$ and $\psi_2$
with independent sets of variables. The first formula $\psi_1$ is
a finite sequence of polynomial equalities with rational
coefficients expressing the terms $u_r(e_i)$ appearing in $j_s
(X^{e_i} \circ h)$, for $1 \leq r \leq s - \langle \nu, e_i
\rangle$, in terms of the variables $\{u_r(e_{k_i}) \mid 1 \leq i
\leq l, \ 1 \leq r \leq s - \langle \nu, e_{k_i} \rangle \}$. We
deduce that $\chi_f ([\psi_1]) = \L^{sl
 -\ord_{\J_{l}}(\nu)}$.
The second formula comes from the effect on the initial
coefficients $c(e_i)$ for $e_i \in \ell_\nu^s$, of the operation
taking the $s$-jet of an arc.  This operation is described  by
taking the image by the map
\[\Psi:T':=\mathrm{Spec} \, k[c(e_i)^{±1}]_{e_i\in\ell_\nu^s}\rightarrow
T:=\mathrm{Spec}\,  k[c(e_i)^{±1}]_{\langle \nu, e_i
\rangle \leq s}, \] of the point determined by $h \in H^*_\nu$.
The map $\Psi$ is the unramified covering of $l$-dimensional
algebraic tori determined by the inclusion $M_\nu^s \subset
\ell_\nu^s \cap M$ of index $q(\nu,s)$ of rank $l(\nu,s)$ lattices. Thus the second
formula is equivalent to $\psi_2 : ( \exists y \in T') (  \Psi (y) = x, \mbox{
and } x \in T)$, hence by Lemma \ref{Lema0} we get that $\chi_f
([\psi_2]) =\frac{1}{q(\nu,s)}(\L-1)^l$. \hfill $\Box$

\section{Sequences of convex piece-wise linear functions and
fans} \label{fan}

Let $\s \subset N_\R$ be a rational convex polyhedral cone of
dimension $d = \dim N_\R$. Consider a sequence of piece-wise
linear continuous functions
\[
h_p: \s \rightarrow \R, \mbox{ for } 1 \leq p \leq m,
\]
such that $h_p (\s \cap N) \subset \Z$, and
\begin{equation} \label{desi}
0 \leq h_1 (\nu) \leq \dots\leq h_{m} (\nu) \quad \forall \nu \in
\s.
\end{equation}
By convenience we set $h_0 (\nu)  =0$ and $h_{m+1} (\nu ) = +
\infty$.  We denote by  $\Xi_0$ the fan consisting on the faces of
$\s$ and  by $\Xi_p$ the coarser fan such that the restriction of
$h_{p}$ to $\eta$ is linear for any cone $\eta \in \Xi_p$ for  $1
\leq p \leq m$. In addition we assume that for any cone $\eta \in
\Xi_{p-1}$ the restriction $h_{p | \eta}$ is {\em upper convex},
that is $h_p (\nu) + h_p (\nu') \leq h_p (\nu + \nu')$ for all
$\nu, \nu' \in \eta$.

\begin{Not} \label{fan-not1}
For $0 \leq p \leq m$  and for $\eta \in \cap_{r =1}^p \Xi_r $ we
set
\[
\eta (\underline{h}, p ) := \{ (\nu, s)\in N_\R\times  \R_{\geq 0}
\mid \nu \in \stackrel{\circ}{\s} \cap \stackrel{\circ}{\eta}, \,
h_p (\nu) \leq s < h_{p+1} (\nu) \}.
\]
\end{Not}

\begin{Lem} \label{fanxi}
The closure  $\bar{\eta} (\underline{h}, p )$ of the set $\eta
(\underline{h}, p )$ is a convex polyhedral cone which is rational
for the lattice $N\times  \Z$.
\end{Lem}
{\em Proof.} If $\eta \in \cap_{r=0}^p \Xi_r$ then the restriction
$h_{j | \eta} \colon \eta \to \R$ is linear if $j = p$ and upper
convex if $j= p+1$. It follows that $\bar\eta(\underline{h},p)$ is
a convex polyhedral cone, rational for the lattice $N \times  \Z$
since $h_p$ and $h_{p+1}$ take integral values on $N$. \hfill $\,
{\Box}$

\begin{Not}  For $0 \leq p \leq m$ and $\eta \in \Xi_p$ we define the following sets: \label{fan-not2}

\begin{enumerate}
\item[(i)] $ A (\underline{h}, p) := \{ (\nu, s)\in N\times  \Z \mid
\nu \in \stackrel{\circ}{\s} ,  \ h_p (\nu) \leq s < h_{p+1} (\nu)
\}$. \item[(ii)] $ A(\underline{h}, p, \eta) := \{ (\nu, s)\in
N\times  \Z \mid \nu \in \stackrel{\circ}{\s}  \cap
\stackrel{\circ}{\eta},  \ h_p (\nu) \leq s < h_{p+1} (\nu) \}$.
\end{enumerate}
\end{Not}

\begin{remark} \label{party}
We have partitions
\[
(\stackrel{\circ}{\s} \cap N) \times  \Z_{\geq 0} = \bigsqcup_{p
=0}^m A (\underline{h}, p) \quad \mbox{ and } \quad
A(\underline{h}, p) =  \bigsqcup_{\eta \in \cap_{r =0}^p \Xi_r }
A(\underline{h}, p, \eta).
\]
\end{remark}

\section{Refinements of partitions}

\label{Combinatoria}

We apply the procedure of Section \ref{fan} to both sequences
$\underline{\f} = (\f_1, \dots, \f_d)$ and $\underline{\varphi} =
(\varphi_1, \dots, \varphi_n)$ (see Lemma \ref{ine}).

\begin{remark}
Notice that the sequence of fans associated to $\underline{\f}$
(resp. $\underline{\varphi}$) is $\cap_{r=0}^i \Sigma_r$, $i =0,
\dots, d$ (resp. $\cap_{r=0}^i \Theta_r$, $i =0, \dots, n$), where for
convenience we denote  by $\Sigma_0$ or by  $\Theta_0$ the fan
consisting of the faces of the cone $\s$.
\end{remark}

\begin{Lem}  \label{lq} If    $A(\underline{\varphi}, j , \theta) \ne \emptyset$
for some $1\leq j \leq n$ and $\theta \in \cap_{r=1}^j \Theta_r$
 (cf. Notation \ref{fan-not2})  then the restriction of the functions $(\stackrel{\circ}{\s} \cap N) \times  \Z_{>0} \to
\Z_{\geq 0}$ given by
\[
(\nu, s) \mapsto l(\nu,s), \quad \mbox{ and } \quad (\nu, s)
\mapsto q(\nu,s),
\]
to the set $A(\underline{\varphi}, j , \theta)$ are constant
functions. We denote their values on the set
 $A(\underline{\varphi}, j , \theta)$ by
 $l(j, \theta)$ and $q(j, \theta)$ respectively.
\end{Lem}
{\em Proof.}
 By elementary properties of Minkowski sums every
vector in the relative interior of $\theta$, for  $\theta \in
\cap_{r=1}^j \Theta_r$, defines the same face $\mathcal{F}_{r,
\theta}$ of the polyhedron $\mathcal{N} (\cc_r)$ for $1\leq r\leq j$. Suppose that
$\nu, \nu' \in \stackrel{\circ}{\theta}$ and  $(p_1, \dots, p_n)$
and $(p_1', \dots, p_n')$ are two permutations of $(1, \dots, n)$
such that the inequalities
\begin{equation} \label{star1}
\langle \nu, e_{p_1} \rangle \leq \cdots \leq \langle \nu, e_{p_n}
\rangle \mbox{ and } \langle \nu', e_{p_1'} \rangle \leq \cdots
\leq \langle \nu', e_{p_n'} \rangle,
\end{equation}
hold. We prove first that
\begin{equation} \label{star2}
 \langle \nu, e_{p_1'} \rangle \leq \cdots
\leq \langle \nu, e_{p_j'} \rangle.
\end{equation}
By definition, for any $1 \leq r \leq j$ we have that $\ord_{\cc_
r} (\nu) = \langle \nu, u_r \rangle$ for any $u_r \in
\mathcal{F}_{r, \theta}$.  We get from Lemma \ref{ine} that the
vectors $u_r := e_{p_1} + \cdots +  e_{p_r}$ and $u_r' := e_{p_1'}
+ \cdots + e_{p_r'}$ belong to $ \mathcal{F}_{r, \theta}$ for $1\leq r\leq j$. This
implies (\ref{star2}).

If $(\nu,s) \in A(\underline{\varphi}, j , \theta)$ then by Lemma
\ref{ine} we obtain that  $\varphi_j (\nu) = \langle \nu, e_{p_j}
\rangle \leq s < \varphi_{j+1} (\nu)$.
We deduce that if $(\nu,s)$ and
$(\nu',s')$ belong to $A(\underline{\varphi}, j , \theta)$ then
\begin{equation} \label{star3}
\{ e_i \mid 1 \leq i \leq n, \ \langle \nu, e_i \rangle \leq s \}
= \{ e_i \mid 1\leq i \leq n, \ \langle \nu', e_i \rangle \leq s'
\} = \{ e_{p_1}, \dots, e_{p_j} \}.
\end{equation}
Since (\ref{star3}) spans the lattice $M_\nu^s$ and the vector
space $\ell_\nu^s$ we get that the sublattices  $\ell_\nu^s \cap M
$ and $M_\nu^s$ are independent of the choice of $(\nu,s)$ in
$A(\underline{\varphi},j,\theta)$. This implies the constancy of
the functions $l$ and $q$ on $A(\underline{\varphi},j,\theta)$.
\hfill $\ {\Box}$

\begin{remark}  \label{compara}
If $1 \leq l \leq d$  and if $\t \in \cap_{r=1}^l \Sigma_r$ we
denoted in \cite{CoGP} the set $A(\underline{\phi}, l )$  (resp.
$A(\underline{\phi}, l, \t )$) by $A_{l}$ (resp. by $A_{l, \t}$).
The map $l(\nu,s)$ is also constant on the sets
of the form $A(\underline\phi,l,\t)$ for $\t\in\cap_{i=0}^l\Sigma_i$, see Lemma 5.7 \cite{CoGP}.
\end{remark}

By Remark \ref{party} we have two partitions:
\begin{equation}
(\stackrel{\circ}{\s}\cap N)\times  \Z_{>0} =
 \bigsqcup_{j=0}^n \bigsqcup_{\theta\in\cap_{r=0}^j\Theta_r} A(\underline{\varphi},j,\theta)
 \ \mbox{ and }\ (\stackrel{\circ}{\s}\cap N)\times  \Z_{>0}= \bigsqcup_{l=0}^d \bigsqcup_{\eta
 \in\cap_{i=0}^l\Sigma_i}
A(\underline{\phi},l,\eta), \label{ast1}
\end{equation}
 associated to the sequences
$\underline{\varphi}$ and $\underline{\phi}$.

\begin{Pro} \label{subdivision}
If $\theta(\underline{\varphi}, j ) \ne \emptyset$ for some $1\leq
j \leq n$ and $\theta \in \cap_{r=1}^j \Theta_r$ then there exists
a unique cone $\t \in  \cap_{r=1}^{l(j, \theta)} \Sigma_r$ such
that $\theta \subset \t$ and
\begin{equation} \label{inc}
A( \underline{\varphi}, j , \theta) \subset  A(
\underline{\f}, l(j, \theta) , \t).
\end{equation}
\end{Pro}
{\em Proof.}
Given $(\nu,s)$ and $(\nu',s')$ in $A( \underline{\varphi}, j
, \theta)\subset(\stackrel\circ\s\cap N)\times\Z_{>0}$, we deduce from (\ref{ast1}) that there exist cones $\t\in\cap_{i=0}^l\Sigma_i$ and $\t'\in\cap_{i=0}^{l'}\Sigma_i$ for integers $0 \leq l, l' \leq
d$ such that $(\nu,s) \in A(\underline{\phi},l,\t)$ and $
(\nu',s') \in A(\underline{\phi},l',\t')$. By Lemma 5.7
\cite{CoGP} we have that  $l = l (\nu,s) $ and $ l' = l (\nu',s')$,
and then $l = l'$ by (\ref{star3}). Notice then that $l=l(j,\theta)$ by definition in Lemma \ref{lq}.

 Let  $(p_1,
\dots, p_n)$ and $(p_1', \dots, p_n')$ be two permutations of $(1,
\dots, n)$ such that (\ref{star1}) holds. Then we can apply the
method given in Proposition 5.1 \cite{CoGP}  to determine the
value of $\ord_{\J_i} (\nu)$, $1 \leq i \leq l(\nu,s)$. Moreover,
it is enough to apply this on the set (\ref{star3}) instead of on
$\{ e_1, \dots, e_n \}$.
We deduce from (\ref{star2}) that $\nu$
and $\nu'$ define the same face of $\mathcal{N}(\J_i)$ for $1 \leq
i \leq l(\nu,s)$. This is equivalent to the equality $\t = \t'$. We have proven (\ref{inc}) and, as a consequence, the inclusion $\theta\subset\t$ holds.
 \hfill $\ {\Box}$

\begin{definition}\label{sim}  (see Definition 8.1  and Remark 8.6 \cite{CoGP}). We consider the equivalence relation
 $\sim$ defined on the set $(\stackrel{\circ}{\s} \cap
N) \times  \Z_{>0}$ by:
\[
(\nu, s) \sim (\nu', s') \quad \Leftrightarrow  \quad s= s', \,
\ell_\nu^s = \ell_{\nu'}^s \mbox{ and } \nu_{|\ell_\nu^s} =
\nu'_{|\ell_{\nu'}^s}.
\]
\end{definition}

\begin{Lem} \label{classes}
The set $A(\underline{\varphi}, j , \theta)$ is union of
equivalence classes by the relation $\sim$ of Definition
\ref{sim},
 for $1\leq j \leq n$
and $\theta \in \cap_{r=1}^j \Theta_r$. Moreover we have that
\begin{equation}
A(\underline{\phi},l,\t)/_\sim=\bigsqcup_{\theta\in\cap_{r=1}^j\Theta_r,\,
l(j,\theta)=l}^{\theta\subset\t}A(\underline{\varphi},j,\theta)/_\sim.
\label{sim2}
\end{equation}
\end{Lem}
{\em Proof.}
By (\ref{ast1}) and Proposition \ref{subdivision} it follows that $A(\underline\phi,l,\t)=\bigsqcup_{\theta\in\cap_{r=1}^j\Theta_r,\,
l(j,\theta)=l}^{\theta\subset\t}A(\underline{\varphi},j,\theta)$. If $(\nu,s)$ belongs to $A( \underline{\varphi}, j ,
\theta)$
 and $(\nu,s) \sim (\nu',s)$  then (\ref{star3}) holds. The vectors $\nu$ and $\nu'$ define the same face of $\mathcal N(\mathcal C_r)$ for $1\leq r\leq j$, and therefore $\nu'\in\mbox{int}\theta$. Since $\varphi_j (\nu') \leq s < \varphi_{j+1} (\nu ')$ we conclude that $(\nu',s)\in A(\underline\varphi,j,\theta)$. \hfill
$\Box$

\section{The structure of the series $\Par^{(Z^\Lambda, 0)} (T)$}
\label{ArithToric}

We consider the following auxiliary Poincaré series:

\begin{equation}
   \Par(\Lambda):=
    \displaystyle\sum_{s\geq 0} \chi_f ( [j_s(H_\Lambda)\setminus\displaystyle\bigcup_{0 \neq\theta\leq\s}
    j_s(H_{\Lambda\cap\theta^\bot}) ]_f ) T^s \in \kmotq
    [[T]].
\label{auxiliarP}
 \end{equation}

 Notice that the Poincaré series $\Par(\Lambda)$ measures the class of the formula defining
 the set of jets of arcs  with origin in $0$  which are not jets of arcs factoring through proper orbit
 closures of the toric variety $Z^\Lambda$.
\begin{Pro} \label{descompPgeom}   We have that
$
\Par^{(Z^\Lambda,0)}(T)=\displaystyle\sum_{\theta\leq\s} \Par
(\Lambda\cap\theta^\bot)$.
\end{Pro}
It follows from Proposition \ref{descompPgeom} that in order to
describe the motivic series $\Par^{(Z^\Lambda,0)}(T)$
it is enough to describe the form of the auxiliary series $\Par
(\Lambda)$ for any semigroup $\Lambda$.

\begin{remark}
In the normal case the equality $j_m(H_\Lambda)=j_m(H_\Lambda^*)$
holds for all $m \geq 0$, see \cite{Nicaise1}, but this property fails
in general. \label{ejemplo}
\end{remark}

We recall the following definition from \cite{CoGP}.

\begin{definition} \label{Sk}
(see \cite{CoGP}  Definition 8.9) If $1 \leq l \leq d$  the set
$\mathcal{D} _l $ is the subset of cones $\t \in \bigcap_{i= 1}^l
\Sigma_i$ such that the face $\mathcal{F}_\t$ of $\mathcal{N}
(\J_l)$ is contained in the interior of $\s^\vee$.
We denote by $|\mathcal{D}_l|$ the set $\cup_{\t \in
\mathcal{D}_l} \t$.
\end{definition}

\begin{Pro} \label{29}
Let us fix an integer $s_0  \geq 1$. The set $j_{s_0}
(H_\Lambda^*)\backslash\bigcup_{0 \neq \theta\leq\s}j_{s_0}
(H_{\Lambda\cap\theta^\bot})$ expresses as a finite disjoint union
of locally closed subsets, as follows:
\begin{equation} \label{trunc-s}
j_{s_0}(H_\Lambda^*)\setminus\displaystyle\bigcup_{ 0
\neq\theta\leq \s}
j_{s_0}(H_{\Lambda\cap\theta^\bot})=
\bigsqcup_{j=1}^n       \,
\bigsqcup_{\theta \in \cap_{r=1}^j \Theta_r}^{\theta \subset |\mathcal{D}_{l(j, \theta)}| }    \,
\bigsqcup_{[ (\nu, s_0) ] \in A(\underline{\varphi}, j, \theta)/_{\sim}}  j_{s_0} (H_{\Lambda,\nu}^*).
\end{equation}
\end{Pro}
\textit{Proof.}
This partition follows from the partition given in Proposition 8.11
\cite{CoGP} by using  Formula (\ref{sim2}) (see Remark \ref{compara}).   \hfill $\Box$

If $s_0 \geq 1$ the coefficient of $T^{s_0}$ in the auxiliary
series $P(\Lambda)$ is obtained by applying the
map $\chi_f$ to the class of the formula defining (\ref{trunc-s}).
Then we determine this class by using  Proposition  \ref{ProArithKey}.

We introduce the following auxiliary series for $\theta \in \cap_{r=1}^j \Theta_r$:
\begin{equation}  \label{Pktau}
P_{\underline\varphi, j, \theta}  (\Lambda): = (\L-1)^{l(j,\theta)}
\,
\sum_{s\geq 1}
\,
\sum_{[ (\nu, s) ]
\,
\in A({\underline\varphi,j, \theta} )/_{\sim}}
\,
\L^{ l(j, \theta)s  -\ord_{\J _{l (j, \theta)}}(\nu)}T^s.
\end{equation}

We deduce the following  Proposition from Proposition \ref{29} and Formula (\ref{Pktau}).
\begin{Pro}    \label{30}
We have that
\begin{equation} \label{28}
\Par (\Lambda) = \sum_{j=1}^n
\      \,
\sum_{\theta \in \cap_{r=1}^j \Theta_r}^{\theta \subset |\mathcal{D}_{l(j, \theta)}| }
\          \,
  \frac{1}{q(j,\theta)} \  P_{\underline\varphi, j, \theta}  (\Lambda).
\end{equation}
\end{Pro}

\section{The rational form of some generating series}
\label{RatGenSeries}

In this Section we fix an integer  $1\leq j \leq
n$ and a cone $\theta \in \cap_{r=1}^j \Theta_r$ such that
$A(\underline{\varphi}, j , \theta) \ne \emptyset$. For simplicity
we denote by $l$ the integer $l( j, \theta)$ defined in Lemma
\ref{lq} and by $\t$ the unique  element of the fan $\cap_{r =1} ^l
\Sigma_r$ such that (\ref{inc}) holds.

Since $\theta \subset \t \subset \cap_{r=1}^l \Sigma_r$ the
restriction of $\f_r$ to $\theta$ is a linear function of the form
\[
(\f_r)_{| \theta } (\nu) = \langle \nu, e_{i_r} \rangle, \mbox{
for } r=1, \dots, l
\]
where $\{ i_1, \dots, i_l \} \subset \{ 1, \dots, n \}$.

Consider the lattice homomorphisms
\[
µ\colon N \times  \Z \longrightarrow  \Z^{l+1}, \quad (\nu,s )
\mapsto (\langle \nu, e_{i_1} \rangle, \dots, \langle \nu, e_{i_l}
\rangle, s)
\]
and
\[
\p = (\p_1, \p_2) \colon \Z^{l+1} \longrightarrow  \Z^2, \quad
(a_1, \dots, a_{l+1}) \mapsto (l a_{l+1} - a_1 - \dots- a_l,
a_{l+1}).
\]
We set $\xi = \p \circ µ$.

\begin{remark} The homomorphisms $\p, µ$  and $\xi$ were also
considered in \cite{CoGP}. Since $\theta$ is contained in $\t$  we
get that
the kernels of $µ$ and $\xi$ intersect the cone $\theta$ only at the origin.
Similarly by Formula (\ref{inc})
 the inclusion
 $\xi( A( \underline{\varphi}, j, \theta ))
\subset \Z^2_{\geq 0} \setminus \{ (0,0) \}$ holds. See \cite{CoGP} Section 9.
\end{remark}

If $j \ne n$ the {\em lower boundary}  of the cone
$\theta$ is the set
$
\partial_- \theta := \{ (\nu, s) \mid \nu
\in {\theta}, \ s = \varphi_j (\nu) \}$. Notice that $\partial_- \theta$ is a cone since $\theta \in \cap_{r=1}^j \Theta_j$
and then
the function $\varphi_j$ is linear on $\theta$.
The {\em upper boundary} is the set   $\partial_+ \theta(\underline{\varphi}, j) := \{ (\nu, s) \mid \nu
\in {\theta}, \ s = \varphi_{j+1} (\nu)\neq\varphi_j(\nu) \}$.
If $j = n$ then $l=d$ and $\varphi_{n+1} (\nu) = + \infty$ and the upper
boundary is the union of cones spanned by $(0,1) \in
N_\R \times  \R$ and the proper faces of
the cone $\partial_- \theta(\underline{\varphi}, j)$.
The edges  of the cone $\theta(\underline{\varphi}, j)$ are
edges of  $ \partial_- \theta(\underline{\varphi}, j) \cup
\partial_+\theta(\underline{\varphi}, j) $.

\begin{Not}
 If $\r \subset \t$ is a one-dimensional cone rational
   for the lattice $N$ we denote by $\nu_\r $ the primitive integral vector on $\r $,
    that is, the generator of the semigroup $\r \cap N$.
 \end{Not}

\begin{remark} \label{edge} The primitive integral vectors for the lattice $N \times  \Z$ on the edges of the cone
 $\theta$ are
\[
\begin{array}{lcl}
(\nu_\r , \varphi_j (\nu_\r )) &  \mbox{ for } &  \r \leq \theta, \,
\dim \r = 1
\end{array}
\]
together with
\[
\left\{
\begin{array}{lcl}
(\nu_\r , \varphi_{j+1} (\nu_\r ))   \ \mbox{ for }  \r \in
\Theta_{j+1},  \ \r \subset {\theta}, \ \dim \r = 1 \ \mbox{ and }\varphi_j(\nu)\neq\varphi_{j+1}(\nu) & \mbox{ if
}  & j \ne n
\\
(0, 1) &  \mbox{ if } & j = n.
\end{array} \right.
\]
Then notice that
\begin{equation}
\xi ( \nu, s) = \left\{
\begin{array}{lcl}
(l \varphi_j (\nu_\r ) - \ord_{\J_l} (\nu_\r ) , \varphi_j (\nu_\r )
) & \mbox{ if } & (\nu, s) = (\nu_\r , \varphi_j (\nu_\r ))
\\
(l \varphi_{j+1} (\nu_\r ) - \ord_{\J_l} (\nu_\r ) , \varphi_{j+1}
(\nu_\r ) )  &  \mbox{ if } & (\nu, s) = (\nu_\r , \varphi_{j+1}
(\nu_\r ))
\\
(d, 1) &  \mbox{ if } & (\nu, s) = (0,1).
\end{array} \right.
\label{xi}
\end{equation}
\end{remark}

\begin{definition}
Suppose that $A(\underline\varphi,j,\theta)\neq\emptyset$. We denote by
$B_{\underline\varphi,j,\theta}(\Lambda)$ the finite set:
\[
 \{(l\varphi_j(\nu_\r )-\ord_{\J_l}(\nu_\r ),\varphi_j(\nu_\r )) \mid {\r \leq\theta},  {\mathrm{dim}\r =1}  \}
\]
\begin{center}$ \cup \left\{ \begin{array}{lcc}
\{(l\varphi_{j+1}(\nu_\r )-\ord_{\J_l}(\nu_\r ),\varphi_{j+1}(\nu_\r )) \mid
\r \in\Theta_{j+1}^{(1)}, \r \subset\theta
 \}& \mbox{ if } & j\neq n,\\

\\

\{(d,1)\} & \mbox{ if } & j=n.\\

\end{array} \right.$\end{center}
\label{ConjB}
\end{definition}

\begin{definition}

If $A \subset \Z^{l+1}$ is a set we denote by $F_A (x) := \sum_{a \in A} x^a$ the
\textit{generating function} of $A$ (see Section 12 of \cite{CoGP})
\label{defFG}.
\end{definition}

\begin{Pro}
We have the following equality:
\begin{equation}       \label{a1}
 P_{\underline{\varphi}, j, \theta} (\Lambda) = (\L-1)^{l(j,\theta)} \sum_{a \in µ(A
(\underline{\varphi}, j , \theta))} \L^{\p_1 (a) } T^{\p_2(a)} \in
\Z [\L][[T]].
\end{equation}
  There exists a polynomial
$R_{\underline{\varphi}, j, \theta} \in \Z [\L, T]$ such that $P_{\underline{\varphi}, j, \theta} (\Lambda)$
has the rational form:
\[
 P_{\underline{\varphi}, j, \theta} (\Lambda) =
R_{\underline{\varphi}, j, \theta} \prod_{(a,b)\in B_{\underline\varphi,j,\theta}(\Lambda)}(1-\L^aT^b)^{-1}.
\]
\label{ratForm}
\end{Pro}
\textit{Proof}.
The map $µ$ defines a bijection
\[
A (\underline{\varphi}, j , \theta) /_{\sim}  \longrightarrow  µ
(A (\underline{\varphi}, j , \theta)), \quad [ (\nu,s) ] \mapsto
µ(\nu,s).
\]
(see  Lemma 9.3 \cite{CoGP} and Lemma \ref{classes}).
Then the equality (\ref{a1})  follows from the definitions.

We denote by  $\pi_*:k[[\Z^{l+1}]]\rightarrow k[[\L,T]]$ the monomial transformation defined by
$\pi_*(x^a):=\L^{\pi_1(a)}T^{\pi_2(a)}$ for $a \in \Z^{l+1}$.
 Then we get that
\[ P_{\underline\varphi,j,\theta}(\Lambda)=   (\L-1)^{l(j,\theta)}
 \pi_* (F_{µ(A_{\underline\varphi,j,\theta})}(x)).\]

We apply the Theorem 12.4 of
\cite{CoGP}.  We obtain that the denominator of a rational form of
$F_{µ(A_{\underline\varphi,j,\theta})}(x)$ consists of products of terms
$1-x^{µ(b)}$ for $b$ running through the primitive integral
vectors in the edges of the closure of the cone $\theta(\underline\varphi,j)$. Then the
result follows by Remark \ref{edge} and Definition \ref{ConjB}. \hfill $\ {\Box}$

\section{Main results}
\label{main}

We summarize the main results of the paper.

\begin{definition}  $\,$    \label{intQ}

\begin{enumerate}
 \item[(i)]  If $0 \leq \eta < \s$ then $\Bar (\Lambda \cap \eta^\bot)$ is the finite subset of $\Z^2_{\geq 0}$
given by Definition    \ref{ConjB}   when we replace $\Lambda$ by the semigroup $\Lambda \cap \eta^\bot$.
If $\eta =\sigma$ we set    $\Bar(\Lambda \cap \s^\bot) : = \{ (0,1)\}$.
We define the finite sets:
\[
  \Bar (\Lambda) := \bigcup_{\theta \in \cap_{r=1}^j \Theta_r, \,  \theta \subset | \mathcal{D}_{l(j, \theta)} | }^{1\leq j \leq n}
B_{\underline{\varphi}, j, \theta}{(\Lambda)}    \quad \mbox{ and } \quad
B_{\mathrm{ar}, \Lambda} : = \bigcup_{0 \leq \eta \leq \s} \Bar (\Lambda \cap \eta^\bot).
\]
\item[(ii)]
We define the integer
\begin{equation}      \label{qL}
      q(\Lambda):=\mathrm{lcm}\{q(j,\theta) \mid
 \theta\in\cap_{r=1}^j\Theta_r,\ \theta\subseteq|D_{l(j,\theta)}|, 1\leq j\leq n\}.
\end{equation}
If $\eta < \s$ then $q(\Lambda \cap \eta^\bot)$ is the number obtained by replacing $\Lambda$ by
$\Lambda \cap \eta^\bot$ in (\ref{qL}).
We set $q(\Lambda\cap\s^\bot):=1$. We define also the integer
\[
 q_\Lambda :=  \mathrm{lcm} \{     q(\Lambda \cap \eta^\bot) \mid \eta \leq \s \}.
\]

\end{enumerate}
\end{definition}

\begin{The}   \label{PLambdaRac}
Suppose that the field $k$ contains all $q(\Lambda)$-th roots of unity.
Then there exists a polynomial  $Q_{\mathrm{ar}} (\Lambda)  \in\Z[\L,T]$ such that
\[
\Par (\Lambda) = \frac{1}{q(\Lambda)}Q_{\mathrm{ar}} ( \Lambda )   \prod_{(a, b) \in \Bar(\Lambda) }
(1-\L^{a}T^{b})^{-1}.   \]
\end{The}
\textit{Proof.} This follows from Propositions \ref{30} and \ref{ratForm}.    \hfill $\Box$

\begin{Not}
 If $\eta < \s$ then the polynomial
$ Q_{\mathrm{ar}} ( \Lambda \cap \eta^\bot )$ is obtained
from Theorem  \ref{PLambdaRac} by replacing $\Lambda$ by the semigroup
$\Lambda \cap \eta^\bot$. We set     $  Q_{\mathrm{ar}} ( \Lambda \cap \s^\bot ) := 1$.
\end{Not}

\begin{Cor}   \label{P-geom}
If the field $k$
contains all $q_\Lambda$-th roots of unity then
there exists a polynomial $Q_{\mathrm{ar}, \Lambda} \in \Z[\L,T]$ such that
such that
\[
 \Par^{(Z^\Lambda, 0)}(T) = \frac{1}{q_\Lambda}  Q_{\mathrm{ar}, \Lambda}   \prod_{(a, b) \in
B_{\mathrm{ar}, \Lambda}} (1-\L^{a}T^{b})^{-1}.
\]
Moreover,
 we have the equality:
\begin{equation}    \label{ratexp}
    \Par^{(Z^\Lambda, 0)}(T) = \sum_{\eta \leq \s}\frac{1}{q(\Lambda\cap\eta^\bot)}   Q_{\mathrm{ar}}
(\Lambda \cap \eta^\bot)
\prod_{(a, b) \in \Bar (\Lambda \cap \eta^\bot) }
(1-\L^{a}T^{b})^{-1}.
\end{equation}
\end{Cor}
\textit{Proof.}  The result follows by     Theorem \ref{PLambdaRac}   and Proposition
\ref{descompPgeom}.\hfill $\Box$

     We can compare at this moment the series $\Pgeom^{(Z, 0)} (T)$ and $\Par^{(Z,0)}(T)$ (see Definition  \ref{par-geom}).
In \cite{CoGP} we introduced the series
\begin{equation}
   \Pgeom(\Lambda):=
    \displaystyle\sum_{s\geq 0} \chi_c ( [j_s(H_\Lambda^*)\setminus\displaystyle\bigcup_{0 \neq\theta\leq\s}
    j_s(H_{\Lambda\cap\theta^\bot}) ] ) T^s \in \kmotq
    [[T]],
\label{auxiliarG}
 \end{equation}
and we proved that
\[
\Pgeom^{(Z^\Lambda,0)}(T)=\displaystyle\sum_{\theta\leq\s} \Pgeom (\Lambda\cap\theta^\bot).
\]
\begin{Pro}   If the field $k$
contains all $q(\Lambda)$-th roots of unity, then  \label{compara2}
\[ \Par (\Lambda)  -  \Pgeom (\Lambda)  =     \sum_{j=1}^n
\      \,
\sum_{\theta \in \cap_{r=1}^j \Theta_r}^{\theta \subset |\mathcal{D}_{l(j, \theta)}| }
\          \,
 ( 1- \frac{1}{q(j,\theta)} ) \  R_{\underline{\varphi}, j, \theta} \ \prod_{(a,b)\in B_{\underline\varphi,j,\theta}(\Lambda)}(1-\L^aT^b)^{-1}.
  \]
\end{Pro}
\textit{Proof.}
This follows from Proposition \ref{ratForm}, Formula (\ref{auxiliarG}),   Theorem
\ref{PLambdaRac}, and the results in \cite{CoGP} for $\Pgeom (\Lambda)$. \hfill $\Box$

\begin{Cor}
If for every integer $1\leq l\leq d$, and any vertex $v$ of the
Newton polyhedra $\mathcal N(\J_l)$ there exists a subset $I_v
\subset \{ 1, \dots, n \}$ of $l$ elements such that $ v =
\sum_{i\in I_v} e_i$ and the vectors
$e_i, i \in I_v$, form part of a basis of $M$ then the series
$\Par^{(Z^\Lambda, 0)}(T)$ and $\Pgeom^{(Z^\Lambda, 0)}(T)$ coincide .
\label{Pg=Pa}
\end{Cor}

{\em Proof.} This condition implies that $q(\nu,s)=1$ for every
$(\nu,s)\in(\stackrel{\circ}{\s} \cap N)\times  \Z_{> 0}$. By Proposition
\ref{compara2}  we get that $ \Par (\Lambda) =  P_{{\mathrm{geom}}}(\Lambda)$. Now
for any face $\eta \leq \s$ the vertices of the Newton polyhedra
of the logarithmic jacobian ideals of $\Lambda \cap \eta^\perp$
are also vertices of the logarithmic jacobian ideals of $\Lambda$.
The hypothesis implies that $\Lambda \cap \eta^\bot$ spans the lattice $M \cap
\eta^\perp$ and then also that $ \Par (\Lambda \cap \eta^\perp) =
P_{{\mathrm{geom}}}(\Lambda \cap \eta^\perp )$. \hfill $\Box$

\begin{remark}
Corollary \ref{Pg=Pa} is a generalization of Nicaise condition in
the case of normal toric varieties (see Theorem 1
\cite{Nicaise1}).
\end{remark}

\begin{remark}     \label{inter}
The coordinates of the vectors in the set  $B_{\underline\varphi,j,\theta}(\Lambda)$   can be described geometrically
in terms of the ideals $\cc_j$ and $\J_l$, for $l = l(j, \theta)$.
Let $\p_j: Z_j \to Z^\Lambda$ be the composite of the normalization of $Z^\Lambda$ with the
toric modification defined by the subdivision $\cap_{r=1}^j \Theta_r $ of $\s$. The modification
$\p_j$ is the minimal toric modification which factors through the normalized blowing up with center
$C_r$, for $r =1, \dots,j$. If $\r$ is an edge of $\theta$ the orbit closure $E_\r$ of the orbit
associated to $\r$ on $Z_j$  has codimension $1$. We denote by $v_\r$ the divisorial valuation defined by $E_\r$.
It verifies that $v_\r (X^m) = \langle \nu_\r, m \rangle$, for $m \in M$. The pull back $\p^*_j (\mathcal{I})$ of a monomial
ideal $\mathcal I$ of $Z^\Lambda$ is a sheaf of monomial ideals on $Z_j$. The ideals $\p^*_j (\cc_r)$, $r=1, \dots, j$
are locally principal on $Z_j$. Then we get the following identities:
\[
\varphi_j (\nu_\r) = v_\r (  \p^*_j (\cc_j) ) - v_\r (    \p^*_j (\cc_{j-1} ) ),
\varphi_{j+1} (\nu_\r) = v_\r (  \p^*_j (\cc_{j+1}) ) - v_\r (    \p^*_j (\cc_{j} ) ),
 \ord_{\J_l} (\nu_\r) = v_\r ( \p_j^* (\J_l ) ).
\]
Compare with the geometrical description of the set of candidate poles of $\Pgeom^{(Z^\Lambda, 0)} (T)$, see
\cite{CoGP}.
\end{remark}

\section{The normal case}
\label{NormalCase}
In the normal case, when the semigroup $\Lambda$ is saturated, i.e., $\Lambda = \s^\vee \cap M$
we describe the motivic arithmetic series in a simpler way by using that
$j^s(H_\Lambda) = j_s(H_\Lambda^*)$ (see \cite{Nicaise1}).

\begin{Not}      $\,$
\begin{enumerate}
\item [(i)] $\mathcal
A=\bigsqcup_{l=1}^d\bigsqcup_{\t\in\cap_{r=1}^l\Sigma_r}A(\underline\phi,l,\t)/_\sim$.

\item [(ii)] For $s_0\geq 0$ we set $\mathcal
A_{s_0}=\{[(\nu,s)]\in\mathcal A\ |\ s=s_0\}$.
\end{enumerate}
\label{DefA}
\end{Not}

\begin{remark}
The set $\mathcal A_s$ is finite (see Remark 8.2 in
\cite{CoGP}).  By (\ref{ast1}) and Lemma \ref{classes}
we deduce that
$\mathcal A=\bigsqcup_{j=1}^n\bigsqcup_{\theta\in\cap_{r=1}^j\Theta_r}A(\underline\varphi,j,\theta)/_\sim$.
\label{RemA}
\end{remark}

\begin{Pro}   Let us fix an integer $s_0 \geq 1$. The set  $j_{s_0}(H^*)$ expresses as a finite disjoint union of
locally closed subsets as $j_{s_0}(H^*)=\bigsqcup_{[(\nu,s)]\in\mathcal A_{s_0}}j_{s_0} (H_\nu^*)$.
We deduce that $\chi_f([j_s(H^*)]_f)=\sum_{[(\nu,s)]\in\mathcal A_s}\chi_f([j_s(H_\nu^*)]_f)$.
\label{Prop2}
\end{Pro}
{\em Proof.} The first assertion follows by applying the method of Proposition 8.11 of \cite{CoGP}. The second assertion
 is a consequence of the first and Proposition  \ref{ProArithKey}. \hfill $\Box$

\begin{The}  \label{ThNormal}  If $Z^\Lambda$ is normal then we have
\[ \Par^{(Z^\Lambda, 0)} = \sum_{j=1}^n \ \sum_{\theta\in\cap_{r=1}^j\Theta_r}     \ \frac{1}{q(j,\theta)}
\ R_{\underline{\varphi}, j, \theta} \ \prod_{(a,b)\in B_{\underline\varphi,j,\theta}(\Lambda)}(1-\L^aT^b)^{-1}.
\]
\end{The}
{\em Proof.} It is consequence of Proposition \ref{Prop2}, Remark \ref{RemA} and Proposition \ref{ratForm}. \hfill $\ {\Box}$

\begin{Cor}
Suppose that  the affine toric
variety $Z^\Lambda$ is normal.
If $\theta \leq \s$ we denote by
$\s_\theta^\vee$ the image of the cone $\s^\vee$ in
$(M_\theta)_\R$, where $M_\theta := M / \theta^\bot \cap M$ and by
$\Lambda (\theta)$ the saturated semigroup $\Lambda (\theta) :=
(\s_\theta^\vee \cap M_\theta ) \times  \Z_{\geq
0}^{\mathrm{codim}\theta}$. With this notation  we have
\[
\Par^{Z^\Lambda}(T)=
 \sum_{\theta\leq\s}(\L-1)^{\mathrm{codim}\theta}
 \Par^{(Z^{\Lambda(\theta)},0)}(T).
\]
\label{P-locales}
\end{Cor}
{\em Proof.} The proof follows by the same arguments as in Corollary 4.11 \cite{CoGP}.\hfill$\Box$

\section{Examples}
\label{ejemplos}
\subsection{The case of monomial curves} Let $\Lambda \subset \Z_{\geq 0}$ be a semigroup
with minimal system of generators $e_1<e_2<\cdots<e_n$ such that  $\mathrm{gcd} \{ e_1,\ldots,e_n \}=1$.
In this case we have that $ \stackrel\circ\s\cap N=\Z_{>0}$.
If  $q_i: =\mathrm{gcd} \{ e_1,\ldots,e_i \}$ then we obtain that:
\begin{equation}  \label{fmla}
         P_{\mathrm{ar}}^{(Z^\Lambda,0)}(T)=\frac{1}{1-T}+
\frac{\L-1}{1-\L T}\left(\frac{1}{q_1}\frac{T^{e_1}}{1-T^{e_1}}+
\sum_{i=2}^n \frac{q_{i-1}-q_i}{q_{i-1}q_i}\frac{\L^{e_i-e_1}T^{e_i}}{1-\L^{e_i-e_1}T^{e_i}}\right).
\end{equation}

This follows from the results of this paper taking the following observations into account:
\begin{itemize}
\item We have the equality $j_s(H)=j_s(H^*)$.
\item
If $\nu,\nu'\in \Z_{>0}$ verify that $ j_s(H_\nu^*)$, $j_s(H_{\nu'}^*)$ $\ne \{ 0 \}$
then the equality  $j_s(H_\nu^*)=j_s(H_{\nu'}^*)$ implies that $\nu=\nu'$.
\item If $\nu\in \Z_{>0}$ verifies that
$\nu e_i\leq s<\nu e_{i+1}$ then  $q(\nu,s)=q_i$.
\end{itemize}

Then, setting $e_{d+1}:=\infty$, we get the following equality which implies (\ref{fmla}):
\[\Par^{(Z^\Lambda,0)}(T)=\frac{1}{1-T}+\sum_{i=1}^n \sum_{\nu=1}^\infty
\sum_{s=\nu e_i}^{\nu e_{i+1}-1}(\L-1)\frac{1}{q_i}\L^{s-\nu e_1}T^s.\]
\begin{remark}            \label{ex}
The inequalities  $q_1 \geq q_2 \geq \cdots \geq q_n=1$ are not always strict. For instance
if $\Lambda$ is generated by $e_1 = 8, e_2 =18, e_3 = 20$ and $ e_4 = 21$ then we get
$q_1= 8$, $q_2= q_3 =2$, $q_4 =1$. It follows from (\ref{fmla}) that the term $1 - \L^{12}T^{20}$ is
not a factor of the denominator of the series    $P_{\mathrm{ar}}^{(Z^\Lambda,0)}(T)$.
Notice that if $\Lambda'$ is the semigroup generated by  $e_1, e_2$ and $ e_4$
then we obtain from       (\ref{fmla})  that       $P_{\mathrm{ar}}^{(Z^\Lambda,0)}(T) =
P_{\mathrm{ar}}^{(Z^{\Lambda'},0)}(T)$  while the semigroups $\Lambda$ and $\Lambda'$ are not
isomorphic.
In contrast with this  behavior, the motivic series  $P_{\mathrm{ar}}^{(C,0)}(T)$ of a plane branch $(C,0)$  determines
the semigroup of the branch $(C,0)$ (see \cite{DL-JAMS}).
\end{remark}

\subsection{An example of non-normal toric surface}
\begin{figure}[h]
\unitlength=1.2mm
\begin{center}
\begin{picture}(150,50)(0,-8)

\linethickness{0.5mm}
\put(5,0){\line(1,0){30}}
\put(5,0){\line(0,1){30}}
\put(50,0){\line(1,0){30}}
\put(50,0){\line(0,1){30}}
\put(95,0){\line(1,0){30}}
\put(95,0){\line(0,1){30}}

\put(12,-5){$\Theta_1$}
\put(60,-5){$\Theta_1\cap\Theta_2$}
\put(100,-5){$\Theta_1\cap\Theta_2\cap\Theta_3$}

\linethickness{0.3mm}
\put(5,0){\line(2,4){13}}
\put(50,0){\line(2,4){13}}
\put(95,0){\line(2,4){12}}

\put(50,0){\line(3,3){20}}
\put(95,0){\line(3,3){20}}

\put(95,0){\line(1,6){5}}

\put(15,28){$\r_1=(2,5)$}
\put(59,28){$\r_1 =(2,5)$}
\put(63,21){$\r_2=(3,5)$}

\put(98,31){$\r_3=(1,6)$}
\put(107,26){$\r_1=(2,5)$}
\put(112,21){$\r_2=(3,5)$}

\put(15,5){$\theta_{11}$}
\put(8,25){$\theta_{12}$}

\put(65,5){$\theta_{21}$}
\put(60,16){$\theta_{22}$}
\put(53,25){$\theta_{23}$}

\put(109,5){$\theta_{31}$}
\put(103,14){$\theta_{32}$}
\put(100,21){$\theta_{33}$}
\put(96,27){$\theta_{34}$}
\end{picture}
\end{center}
 \caption{The subdivisions $\Theta_1$, $\Theta_1\cap\Theta_2$ and $\Theta_1\cap\Theta_2\cap\Theta_3$
\label{figuraej2}}
\end{figure}
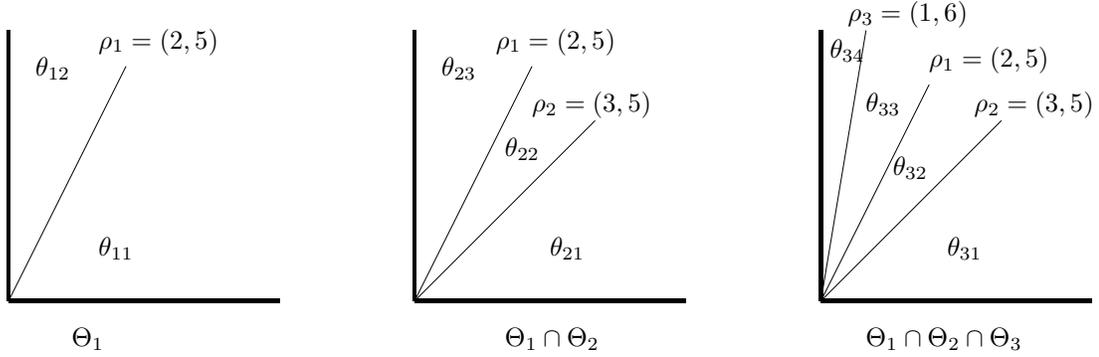

Consider the semigroup $\Lambda$
generated by the vectors $e_1=(5,0), e_2=(0,2), e_3=(0,3)$ and $e_4=(6,2)$.
The cone $\s$ is $\R^2_{\geq 0}$ and the lattice $M$ is equal to $\Z^2$.
We have the semigroups
$\Lambda\cap\eta_1^\bot=(5,0)\Z_{>0}$, and  $\Lambda\cap\eta_2^\bot=(0,2)\Z_{>0}+(0,3)\Z_{>0}$,
 where $\eta_1$ and $\eta_2$ are the one-dimensional faces of $\s$.
By the case of monomial curves we get that:
\begin{center}       $\Par(\Lambda\cap\eta_1^\bot)=\frac{\L-1}{1-\L T}\frac{T}{1-T}$ and
$\Par(\Lambda\cap\eta_2^\bot)=\frac{\L-1}{2(1-\L T)}\left(\frac{T^2}{1-T^2}+\frac{\L T^3}{1-\L T^3}\right)$.
\end{center}
The subdivisions associated with the ideals $\cc_r$, $r=1,\dots, 4$ are indicated in
Figure \ref{figuraej2}

\vspace{-1cm}

In the following table we give the different values of $q(j,\theta)$ and
$l(j,\theta)$,  for $\theta$ in the subdivisions of Figure \ref{figuraej2} and $j$ such that
     $A(\underline\varphi,j,\theta) \ne \emptyset$.   We exclude from this table the cones in
$\theta \in  \cap_{r=1}^4 \Theta_r$ for $j = 4$ since in this case
 $q(4,\theta)=1$ and
$l(4,\theta)=2$.

 \vspace{2mm}
\centerline{\begin{tabular}{|c|c|c|c|c|c|c|}
  \hline
  $j=1$  & $q(1,\theta_{11})=2$ & $q(1,\theta_{12})=5$ &  &  &  & \\
  \hline
   & $l(1,\theta_{11})=1$ & $l(1,\theta_{12})=1$ &  & & & \\
   \hline
  $j=2$ & $q(2,\theta_{21})=1$ & $q(2,\theta_{22})=10$ & $q(2,\theta_{23})=10$ & $q(2,\r_1)=10$ &  &  \\
  \hline
  & $l(2,\theta_{21})=1$ & $l(2,\theta_{22})=2$ & $l(2,\theta_{23})=2$ & $l(2,\r_1)=2$ &  & \\
  \hline
  $j=3$ & $q(3,\theta_{31})=5$ & $q(3,\theta_{32})=5$ & $q(3,\theta_{33})=5$ &  $q(3,\theta_{34})=2$ & $q(3,\r_1)=5$ & $q(3,\r_2)=5$ \\
  \hline
 & $l(3,\theta_{31})=2$ & $l(3,\theta_{32})=5$ & $l(3,\theta_{33})=2$ &  $l(3,\theta_{34})=2$ & $l(3,\r_1)=2$
& $l(3,\r_2)=2$ \\
      \hline
\end{tabular}}

 \vspace{3mm}

Notice that we have $A(\underline\varphi,1,\r_1)=A(\underline\varphi,2,\r_2)=A(\underline\varphi,3,\r_3)=\emptyset$.
In the following table we have filled in the cases  corresponding to the pairs $(a, b) \in  B_{\mathrm{ar}}(\Lambda)$, $(a, b) \ne (2,1)$
in terms of the rays appearing  in the subdivisions of
 Figure \ref{figuraej2}:
 \vspace{2mm}

\centerline{\begin{tabular}{|c|c|c|c|c|c|}
  \hline
 $(a, b) \in  B_{\mathrm{ar}}(\Lambda) $   & $ \nu_{\r_1} = (2,5)$ & $ \nu_{\r_2} = (3,5)$ & $ \nu_{\r_3}= (1,6)$ & $\nu_{\s^\vee\cap\eta_1^\bot} = (1,0)$ &
$\nu_{\s^\vee \cap\eta_2^\bot} = (0,1)$\\
  \hline
  $(2\varphi_2-\ord_{\J_2},\varphi_2)$ & (0,10) & (5,15) &   &   & (2,2)\\
  \hline
  $(2\varphi_3-\ord_{\J_2},\varphi_3)$ & (10,15) & (5,15) &  (19,18) & (5,5) & (2,2) \\
  \hline
  $(2\varphi_4-\ord_{\J_2},\varphi_4)$ & (24,22) & (31,28) & (19,18) & (7,6) & (4,3) \\
  \hline
\end{tabular}}

\vspace{3mm}

It follows that $B_{\mathrm{ar}, \Lambda} =  B_{\mathrm{ar}}(\Lambda) \cup \{ (1,3), (0, 2), (1,1), (0,1) \}$.
We have computed the sum of the series
$\Par^{(Z^\Lambda,0)}(T)$ with the methods of \cite{CoGP}.
We have obtained an irredundant representation of the form $\Par^{(Z^\Lambda,0)}(T)
=R(\L,T)\prod_{(a,b)\in B}(1-\L^aT^b)^{-1}$ with $R(\L,T)\in\Q[\L,T]$
and where $B= B_{\mathrm{ar}, \Lambda} \setminus \{  (24, 22), (31,  28) \}$.

     \bibliographystyle{amsplain}
\def\cprime{$'$}
\providecommand{\bysame}{\leavevmode\hbox to3em{\hrulefill}\thinspace}
\providecommand{\MR}{\relax\ifhmode\unskip\space\fi MR }
\providecommand{\MRhref}[2]{%
  \href{http://www.ams.org/mathscinet-getitem?mr=#1}{#2}
}
\providecommand{\href}[2]{#2}

\end{document}